\documentclass[10pt]{amsart}
\usepackage{amsmath,amsfonts,amsthm,bbm, amssymb}
\usepackage[cmtip, all]{xy}
\usepackage[left=3cm, right=3cm, top = 2.5cm, bottom = 2.5cm ]{geometry}
\usepackage{comment}
 
\usepackage{enumerate}
\usepackage[usenames]{xcolor}
\usepackage{comment}
\usepackage{graphicx}
\usepackage{amssymb}
\usepackage{stmaryrd}
\usepackage{dsfont, mathrsfs,bold-extra}
\usepackage{setspace}

\usepackage[color, notref, notcite, final]{showkeys}     
\definecolor{refkey}{gray}{.5}   
\definecolor{labelkey}{gray}{.5} 
\definecolor{Red}{rgb}{1,0,0}

\usepackage{soul}

\newtheorem{thm}{Theorem}[section]
\newtheorem{prop}[thm]{Proposition}
\newtheorem{lem}[thm]{Lemma} 
\newtheorem{cor}[thm]{Corollary}

\newtheorem{ques}[thm]{Question}

\theoremstyle{definition}

\newtheorem{df}[thm]{Definition}
\theoremstyle{remark}
\newtheorem{rmk}[thm]{Remark}
\newtheorem{remks}[thm]{Remarks}

\newtheorem{exm}[thm]{Example}
\newtheorem{exms}[thm]{Examples}

\numberwithin{equation}{section}

{\hfill$\square$\end{defn}}
{\hfill$\square$\end{remk}}
{\hfill$\square$\end{remks}}
{\hfill$\square$\end{exm}}
{\hfill$\square$\end{exms}}

	\newcommand{\Z}{\mathbb{Z}}

	\renewcommand{\P}{\mathbb{P}}
	\newcommand{\A}{\mathbb{A}}

	\newcommand{\Set}{\mathbf{Set}}

	\newcommand{\Sm}{\mathbf{Sm}}

    \newcommand{\Shv}{\mathbf{Shv}}
    
    \newcommand{\Psh}{\mathbf{Psh}}
    \newcommand{\sPsh}{s\mathbf{Psh}}
    
    \newcommand{\PSh}{\mathbf{Psh}}

\providecommand{\frac}[1]{\operatorname{Frac}(#1)}  
\providecommand{\Spec}[1]{\operatorname{Spec}(#1)} 	%
\newcommand{\Sing}{\operatorname{Sing}}              
  
\newcommand{\id}{\operatorname{id}}   			
\renewcommand{\hom}{\operatorname{Hom}}				 
\renewcommand{\H}{\operatorname{H}}				     
\newcommand{\colim}{\operatorname{colim}}		 
\renewcommand{\lim}{\operatorname{lim}}			 

	\newcommand{\tensor}{\otimes}
	\renewcommand{\cong}{\simeq}
	
	\newcommand{\ol}{\overline}
	\renewcommand{\ul}{\underline}

\newcommand{\cC}{\mathcal{C}}
\newcommand{\cD}{\mathcal{D}}

\newcommand{\cH}{\mathcal{H}}
\newcommand{\cI}{\mathcal{I}}

\newcommand{\cL}{\mathcal{L}}
\newcommand{\cM}{\mathcal{M}}

\newcommand{\cO}{\mathcal{O}}

\newcommand{\cS}{\mathcal{S}}

\renewcommand{\epsilon}{\varepsilon}
\renewcommand{\phi}{\varphi}

\newcommand{\kX}{\mathscr{X}}
\newcommand{\kXc}{\mathscr{X}^{\rm c}}
\newcommand{\kY}{\mathscr{Y}}
\newcommand{\kZ}{\mathscr{Z}}

\newcommand{\MM}{{\mathbf{\ol{M}}}\mathcal{M}(k)}
\newcommand{\MMb}{{\mathbf{\ol{M}}}\mathcal{M}_\bullet(k)}
\newcommand{\MH}{\mathbf{\ol{M}} \mathcal{H}(k)}
\newcommand{\MHp}{\mathbf{\ol{M}} \mathcal{H}_\bullet(k)}
\newcommand{\Mlog}{\cM_{\rm log}(k)}
\newcommand{\BMlog}{B\cM_{\rm log}(k)}

\newcommand{\Daytensor}{\tensor^{\rm Day}}

\def\Smlog{\mathbf{Sm}_{\rm log}(k)}
\def\BSmlog{\mathbf{BSm}_{\rm log}(k)}
\def\MSmlog{\mathbf{\ol{M}Sm}_{\rm log}(k)}

\def\Sm{\mathbf{Sm}(k)}
\def\pt{\operatorname{pt}}
\def\Map{\operatorname{Map}}
\def\Hom{\mathscr{H}\!\textit{om}}

\def\boxq{{\ol{\square}}}
\newcommand\Itensor{\mathop{\mbox{$I$-$\tensor$}}}
\newcommand\cItensor{\mathop{\mbox{$\mathcal{I}$-$\tensor$}}}
\newcommand{\boxtensor}{\mathop{\mbox{$\ol{\square}$-$\tensor$}}}
\newcommand{\SingI}{\operatorname{Sing}^{\tensor}_{\mathcal{I}}}

\makeatletter
\newcommand{\colim@}[2]{%
  \vtop{\m@th\ialign{##\cr
    \hfil$#1\operator@font hocolim$\hfil\cr
    \noalign{\nointerlineskip\kern1.5\ex@}#2\cr
    \noalign{\nointerlineskip\kern-\ex@}\cr}}%
}
\newcommand{\hocolim}{%
  \mathop{\mathpalette\colim@{\rightarrowfill@\textstyle}}\nmlimits@
}
\makeatother

\makeatletter
\newcommand{\ncolim@}[2]{%
  \vtop{\m@th\ialign{##\cr
    \hfil$#1\operator@font colim$\hfil\cr
    \noalign{\nointerlineskip\kern1.5\ex@}#2\cr
    \noalign{\nointerlineskip\kern-\ex@}\cr}}%
}
\newcommand{\varcolim}{%
  \mathop{\mathpalette\ncolim@{\rightarrowfill@\textstyle}}\nmlimits@
}
\makeatother

\makeatletter
\newcommand{\varvarcolim}{%
  \mathop{\mathpalette\ncolim@{}}\nmlimits@
}
\makeatother

\makeatletter
\newcommand{\holim@}[2]{%
  \vtop{\m@th\ialign{##\cr
    \hfil$#1\operator@font holim$\hfil\cr
    \noalign{\nointerlineskip\kern1.5\ex@}#2\cr
    \noalign{\nointerlineskip\kern-\ex@}\cr}}%
}
\newcommand{\holim}{%
  \mathop{\mathpalette\holim@{\leftarrowfill@\textstyle}}\nmlimits@
}
\makeatother

\def\hocof{\operatorname{hocof}}
\def\hofib{\operatorname{hofib}}
\definecolor{winered}{rgb}{0.8,0,0}
\usepackage[pdfauthor={Federico Binda}, %
				pdftitle={Hmodulus} ,colorlinks=true]{hyperref}
\hypersetup{
	bookmarksnumbered=true,
	linkcolor={winered},
	citecolor=blue,
}

\usepackage{float}
\newcounter{elno}   
\newenvironment{romanlist}{
                         \begin{list}{\roman{elno})
                                     }{\usecounter{elno}}
                      }{
                         \end{list}}
                         
\newcounter{elno-abc}   

\begin{document}

\author{Federico Binda}
\title{A  motivic homotopy theory without $\mathbb{A}^1$-invariance}

\AtEndDocument{\bigskip{\footnotesize%
  (F.~Binda) \textsc{Dipartimento di Matematica ``Federigo Enriques'', Universit\`a degli Studi di Milano, Via C.~Saldini 50, 20133 Milano, Italy} 
  \textit{E-mail address},  \texttt{federico.binda@unimi.it}}}
\thanks{Partially supported by the CRC SFB 1085 ``Higher Invariants'' (University of Regensburg) of the Deutsche Forschungsgemeinschaft during the preparation of the manuscript.}
\maketitle


\begin{abstract}
In this paper, we continue the program initiated by Kahn-Saito-Yamazaki by constructing  and studying an unstable motivic homotopy category with modulus $\ol{\bf M}\cH(k)$, extending the Morel-Voevodsky construction from smooth schemes over a field $k$ to certain diagrams of schemes. We present this category as a candidate environment for studying representability problems for non $\A^1$-invariant generalized cohomology theories.
\end{abstract}
\setcounter{tocdepth}{1}
\tableofcontents

\section{Introduction}
\subsection{}In this paper we start developing a machinery in the spirit of Morel-Voevodsky's homotopy theory of schemes as described in \cite{MV}. Their construction, which is now classical, starts from the fundamental observation that $\mathbb{A}^1$, the ``affine line'',  should play in algebraic geometry the role of the unit interval $[0,1]\subset \mathbb{R}$ in topology. The homotopy invariance property (in the usual sense) of an ordinary cohomology theory gets  translated in the algebraic world in the invariance with respect to $\mathbb{A}^1$. This allows to capture many important invariants, but not all of them. In the present work, we consider a notion of ``invariance'' with respect to a different object.   


The lighthouse that guides our construction is the behavior of $K$-theory for possibly singular schemes. If the homotopy invariance property, i.e.~the existence of an equivalence $\mathbf{K}(S\times \mathbb{A}^1)\cong \mathbf{K}(S)$ (of, say, $S^1$-spectra) does not hold for an arbitrary scheme $S$, the projective bundle formula, as formulated e.g.~in \cite{TT}, is still available with only mild finiteness assumptions. In fact, for any scheme $S$ quasi-compact and quasi-separated, the $K$-theory spectrum $\mathbf{K}(S\times \P^1)$ decomposes, as module over $\mathbf{K}(S)$, as follows
\[\mathbf{K}(S\times \P^1) \cong  \mathbf{K}(S)[H]/(H^2), \quad H = [\cO(-1)], H^0 = [\cO] \in K_0(\P^1_\mathbb{Z}), \]
where the summand $[H] \mathbf{K}(S)$ is supported on the hyperplane at infinity of $S\times\P^1$. The pullback along the projection $\pi\colon S\times \P^1 \to S$ induces then an isomorphism (in the homotopy category of spectra)
\[\pi^* \colon \mathbf{K}(S)\xrightarrow{\simeq}\hocof( \mathbf{K}(S) \xrightarrow{\iota_{\infty, *}} \mathbf{K}(S\times \P^1)). \]
The homotopy cofiber computes the $K$-theory of $S\times \P^1$ modulo the $K$-theory of  $S\cong S\times \{\infty\}$.  If we think of $\P^1$ as being the compactification of $\A^1$, with the point at infinity as boundary, we can imagine that such cofiber works as a replacement for the $K$-theory of the open complement $S\times \mathbb{A}^1$ (in fact, it \textit{is} equivalent to the $K$-theory of the complement when $S$ is regular).
In a way, the lost homotopy invariance is  found again in a form of invariance with respect to the pair $(\P^1, \infty)$, this time without regularity assumptions. 

It looked therefore reasonable to conceive a theory where the basic objects are not \textit{schemes} but rather \textit{pairs} \textit{of schemes}, 
and where the role played by $\A^1$ in Voevodsky's theory is played by the \textit{closed box} or the \textit{cube} \[\ol{\square}=(\P^1, \infty).\]
This is, with some simplifications, the point of view adopted by Kahn-Saito-Yamazaki in \cite{KSY} (although with a different starting point), following original insights of Kerz-Saito. Let us fix a  field $k$. Objects in the category $\underline{\mathbf{M}}\mathbf{Sm}(k)$ of modulus pairs (see \cite[Definition 1.1]{KSY} or Definition \ref{def:mod-pair-KSY}) are ``partial compactifications'' $\ol{X}$ (with the total space not necessarily proper over $k$) of a smooth and separated $k$-scheme $X$ and with the specified datum of an effective, possibly non reduced, divisor ``at infinity'' $X^\infty$ with  support on $\ol{X}\setminus X$. The morphisms are given by certain finite correspondences between two compactifications $\ol{X}$ and $\ol{Y}$ satisfying suitable admissibility and properness conditions, and restricting to actual morphisms of schemes between the open complements $X$ and $Y$.

\subsection{}
 The idea of extending the motivic framework to the setting of pairs of schemes, in order to capture some non $\mathbb{A}^1$-invariant phenomena with the aid of non-reduced Cartier divisors, has been explored in the work of many authors, starting from the pioneering work of Bloch-Esnault on additive Chow groups in the early 2000's (see  \cite{BEAdditive} and \cite{BEAdditiveChow}), continued by Park, Krishna and Levine  \cite{KL}, \cite{KP}, \cite{KP3}, \cite{KP4} and then further generalized in \cite{BS} where the current formulation of higher Chow groups with modulus and of relative motivic cohomology was introduced. 

In \cite{BS}, the word \textit{pair} was used to designate the datum $(\ol{X},D)$ consisting of a smooth and separated $k$-variety and an effective Cartier divisor $D$ on it. The relative motivic cohomology groups $\H^*_\cM(\ol{X}|D, \mathbf{Q}(*))$ are conjectured to describe rationally the relative $K$-groups $K_*(\ol{X},D)$, defined as homotopy groups of the homotopy \textit{fiber} 
\[\mathbf{K}(\ol{X},D) = \hofib(\mathbf{K}(\ol{X}) \to \mathbf{K}(D)) \]
(along the inclusion $D\hookrightarrow \ol{X}$),  and are related to other invariants, such as relative Deligne cohomology and relative deRham cohomology: see \cite[Sections 6 and 7]{BS}. The functoriality of these invariants is for morphisms of pairs in the topological sense, i.e., morphisms $f\colon \ol{X}\to \ol{Y}$ that restrict to morphisms of the closed subschemes $f_D\colon  D_X\to D_Y$.  From this point of view, the cofiber construction sketched above does not fit completely well with the general picture: it is in fact contrasting with the idea of compactification with boundary divisor that was coming from the analysis of $K$-theory of singular schemes.

We are then in front of two forces pulling in opposite directions. On one side, in the hope of finding an easy replacement for homotopy invariance, we would like to have a divisor $\partial X$  on a scheme $\ol{X}$ representing the boundary of an abstract compactification, such as $(\P^1, \infty)$. On the other hand, we would like to have a relative theory for a divisor $D$, that we think as effective subscheme of $\ol{X}$, to which we attach relative invariants such as relative $K$-theory, Chow groups with modulus and other generalized relative cohomologies. All phenomena that are arising from  fiber constructions, rather then  cofiber constructions, and that have different functoriality constraints. Our solution is to incorporate both aspects in one category. 

Instead of working with the category $\mathbf{Sm}(k)$ of smooth and separated schemes over a field $k$ (Voevodsky's model) or with the category of smooth modulus pairs $\underline{\mathbf{M}}\mathbf{Sm}(k)$ over $k$ built out from the insights of Kahn-Saito-Yamazaki, we introduce a category of \textit{modulus data}, $\MSmlog$. Objects of $\MSmlog$ are triples \[M = (\ol{M}, \partial M, D_M),\] where $\ol{M}$ is a smooth and separated $k$-schemes and $\partial M$ and $D_M$ are effective Cartier divisors on $\ol{M}$. The different roles are reflected in the notation that we have chosen. The divisor $\partial M$ is assumed to be a strict normal crossing divisor on $\ol{M}$, that we think of as a log-compactification of $\ol{M}\setminus M$. Insisting on our interpretation of the divisor at infinity as a boundary (and unlike \cite{KSY}), we assume it to be always reduced. The non-reduced piece of information will then come from the second divisor $D_M$. A simpler category of schemes with compactifications, $\Smlog$, is identified with the full subcategory of $\MSmlog$ of modulus data with empty ``modulus divisor'' $D_M$,
\[ u\colon \Smlog \hookrightarrow \MSmlog.\]
Both categories come equipped with a symmetric monoidal structure, modelled on the tensor structure of $\underline{\mathbf{M}}\mathbf{Sm}(k)$ of \cite{KSY}, and there are natural cd-structures inherited from the standard Zariski and Nisnevich cd-structures on the underlying categories of schemes over $k$ (see Section  \ref{subsec:topology-Nis-Zar-modulus-cat}).
 
Although much of our construction can be carried out without this assumption, we will furthermore ask that $\partial M +|D_M|_{\rm red}$ forms a strict normal crossing divisor on $\ol{M}$. This restriction turns out to be useful in the definition of a $K$-theory realization of a modulus datum (as discussed in \cite{BThesis} and in \cite{BMilnorChowMod}), and is consistent with the assumptions of \cite{BS}.
 
Morphisms in $\MSmlog$ are modelled on the two constraints of the ``compactification divisor'' and of the ``relative divisor'' that we mentioned before. We refer the reader to Definition  \ref{def:modulus-data} for details, but we remark here that we follow the opposite convention of \cite{KSY} (whence the overline notation instead of the underline notation). Out of $\MSmlog$ we build, in parallel with the Morel-Voevodsky construction of (unstable) motivic homotopy categories in \cite{MV}, a homotopy category with ``modulus''. We briefly sketch  the main difficulties that we have encountered.  

\subsection{}There is a first evident difficulty in trying to extend Voevodsky's formalism of sites with interval to our context. The multiplication map \[\mu\colon \A^1\times \A^1\to \A^1\] does not extend to a morphism $\P^1\times \P^1\to \P^1$, but only to a rational map. It is, however, defined as correspondence, and as such satisfies a suitable admissibility condition.  
This is the path followed by Kahn-Saito-Yamazaki in \cite{KSY}, that eventually leads to the definition of the (homological) triangulated category $\underline{\mathbf{M}}\mathbf{DM}^{\rm eff}(k, \Z)$, built as suitable localization of the derived category of complexes of presheaves on $\underline{\mathbf{M}}\mathbf{Cor}(k)$, the category of \textit{modulus correspondences}, an enlargement of the category $\underline{\mathbf{M}}\mathbf{Sm}(k)$ introduced above.

If considering presheaves on categories of correspondences is the basic input in the construction of ``derived '' categories of motives, the intrinsically \textit{linear} nature of them does not fit well with the desire of putting homotopy theory in the picture. We  then abandon the idea of extending directly the set of admissible morphisms and we turn to a different approach.

The key observation is the following. There is no direct way of making the closed box $\ol{\square}$ into an interval object in the category $\sPsh(\MSmlog)$ of simplicial presheaves on $\MSmlog$, nor in the category $\sPsh(\Smlog)$ of simplicial presheaves over $\Smlog$. Instead, we consider an auxiliary category $\mathbf{BSm}_{\rm log}(k)$, built as localization of  $\Smlog$ with respect to a suitable class of admissible blow-ups. For making sense of this, we need our base field $k$ to admit strong resolution of singularities (see Section  \ref{subsec:inv-rat-maps} for a precise definition).

In $\mathbf{BSm}_{\rm log}(k)$, the multiplication map $\mu\colon \ol{\square}\tensor \ol{\square}\to \ol{\square}$ is an acceptable morphism, and $\ol{\square}$ is naturally an interval object. Starting from it, we build a non-representable simplicial presheaf $\cI$ in $\sPsh(\MSmlog)$, that is a   $\tensor$-interval object with respect to the convolution product of presheaves (see Sections  \ref{subsec:monoidal-structure-MSmlog}- \ref{sec:def-I-object-sec1} and  \ref{sec:Mot-spaces-and-intervals}). It comes naturally equipped with a  map
\[\ol{\square} \xrightarrow{\eta} \cI,\]
and we define the \textit{motivic model structure with modulus} $\MM$ to be the (left) Bousfield localization of the standard (Nisnevich-local) model structure on $\sPsh(\MSmlog)$ (the category of motivic spaces with modulus) to the class of maps $\kX\tensor \cI\to \kX$. The resulting homotopy category $\ol{\mathbf{M}} \mathcal{H}(k)$ is called the \textit{unstable (unpointed) motivic homotopy category with modulus} (see  \ref{sec:Nis-BG}).

Even though the map $\eta$ is not a weak equivalence (and does not become such after localization), every (Nisnevich-local fibrant) $\ol{\square}$-invariant simplicial presheaf $\kX$ is $\cI$-local. Thus,   the category  $\ol{\mathbf{M}} \mathcal{H}(k)$  can serve the purpose of   representing $\ol{\square}$-invariant theories equally well. We investigate the precise relationship between the $\ol{\square}$-theory, the $\cI$-theory and the $\mathbb{A}^1$-theory in  \ref{sec:Mot-spaces-and-intervals} and in  \ref{sec:comparison-intervals}.

\subsection{}In constructing the $\cI$-model structure on $\MM$, we go through a certain amount of technical work for generalizing Voevodsky's homotopy theory for a site with interval. We generalize his results in two directions. First, we deal with a (closed) monoidal structure on simplicial presheaves that is not the Cartesian product but is induced via Day convolution by the monoidal structure on $\MSmlog$. Second, the interval object we consider is not, as remarked above, representable. Nevertheless, we draw inspiration from the work of Morel-Voevodsky \cite{MV} and of Jardine \cite{JardineSpectra} to construct a replacement for the $\A^1$-singular functor $\Sing^{\A^1}$, that we denote 
\[ \Sing_{\cI}^\tensor(-)\colon \MM\to \MM,\]
and that allows us to reprove in this generality the Properness Theorem \cite[Theorem 3.2]{MV} --- which is one of the  non-trivial structural properties of the Morel-Voevodsky construction.
\begin{thm}[see Theorem  \ref{thm:properness} and  \ref{sec:review-properties-MM}] The $\cI$-$\tensor$-localization $\MM = \MM_{\rm inj}^{\cI-loc}$ of the category of motivic spaces with modulus (over $k$) equipped with the local injective (for the Nisnevich topology) model structure is a proper cellular simplicial monoidal (for the Cartesian product) model category.
    \end{thm}
The proofs are necessarily technical, and cover Sections \ref{sec:sing-functor} and \ref{sec:properness}. After this, we conclude with a characterization of simplicial presheaves satisfying the Brown-Gersten (B.G.) property in our context (\ref{sec:Nis-BG}) and with the following  representability result.

\begin{thm}[see Theorem \ref{thm:formal-representability}] Let $\kX$ be a pointed motivic space with modulus that is Nisnevich excisive  (Proposition  \ref{def:Nisn-excisive-prop}) and $\ol{\square}$-$\tensor$-invariant (in particular, $\cI$-$\tensor$-invariant). Then for any $n\geq 0$ and any modulus datum $M$, we have a natural isomorphism \label{thm:intro-repre}
         \[\pi_n(\kX(M)) \cong [S^n\wedge (M)_+,\kX]_{\MHp}.\]
\end{thm}

We conclude this introduction by mentioning a few possible future developments of the theory. In \cite{BThesis}, we constructed for each modulus datum $M$ a $K$-theory spectrum $\mathbf{K}(M)$, that is automatically $\ol{\square}$-invariant. Unfortunately, the assignment $M\mapsto \mathbf{K}(M)$ is not strictly functorial: if $f\colon M\to N$ is a morphism, we defined a pull-back map $f^*\colon \mathbf{K}(N)\to \mathbf{K}(M)$, but the construction is not canonical (i.e. $(g\circ f)^*$ is not strictly equal to $f^*\circ g^*$ as morphisms of spectra), due to the existence of certain diagrams of maps which do commute only up to homotopy. A suitable refinement of this construction would give rise to a presheaf of spectra on the category of modulus data $\MSmlog$, which will then be automatically an object in the associated category of $S^1$-spectra. An obvious variant of Theorem \ref{thm:intro-repre} could be then applied directly to get the representability of relative $K$-theory in (the $S^1$-stabilization of) $\MHp$. In a different direction, it is important to understand  what is the analogue of the  Morel-Voevodsky localization
property in this context: the most naive guess is probably doomed  to failure, given for example the difficulty of formulating the localization property in the context of higher Chow groups with modulus (see \cite{KrishnaOnCycles} and \cite{RulSaito} for counterexamples, and \cite{BK} for related questions). For reciprocity sheaves (so, in the context of sheaves with transfers), an answer to this question is discussed in \cite{BRS}.

Finally, it is not clear what is the analogue of $\mathbb{P}^1$-spectra (or, in other words, what are the ``spheres'' in this world). This last issue is related to a similar problem in the work of Kahn-Saito-Yamazaki, where only \textit{effective} motivic complexes are considered, and a  duality is lurking beyond the curtains.  We plan to (partially) address these issues in a future work. 

\medskip
{\bf Notations and conventions.} Throughout this paper we fix a base field $k$ for which we assume to have resolution of singularities (see Section \ref{subsec:inv-rat-maps}). Unless specified otherwise, all schemes will be assumed to be separated and of finite type over $k$. We write $\mathbf{Sm}(k)$ for the category of smooth quasi-projective $k$-schemes.
\section{Categories of schemes with moduli conditions}\label{sec:Schemes-with-moduli-conditions}

\subsection{Schemes with compactifications}Let $X$ be a smooth $k$-scheme and let $\partial X$ be a reduced codimension $1$ closed subscheme of $X$ with irreducible components $\partial X_1,\ldots, \partial X_N$. We say that $\partial X$ is a strict normal crossing divisor if for every non empty subset $A$ of $\{1, \ldots, N\}$, the subscheme $\partial X_A = \bigcap_{i\in A} \partial X_i$ is smooth over $k$ and of pure codimension $|A|$ in $X$. We will denote by $\partial X_*$ the set of irreducible components of a strict normal crossing divisor $\partial X$ and we write $|\partial X| = \cup_{i=1}^N \partial X_i$ for the support of $\partial X$. If $T_1,\ldots, T_r$ are smooth integral codimension one subschemes of $X$ such that their union is a strict normal crossing divisor, we say that the set $T_1,\ldots, T_r$ form a normal crossing divisor on $X$. The empty scheme is considered a strict normal crossing divisor of every smooth $k$-scheme. 

\begin{df}\label{def:Smlog}The category $\Smlog$ is the category of pairs $(X, \partial X)$, where $X$ is a smooth $k$-scheme and $\partial X$ is a strict normal crossing divisor on $X$ (possibly empty). A morphism \[f\colon (X, \partial X)\to (Y, \partial Y)\] of pairs in $\Smlog$ is a $k$-morphism $f\colon X\to Y$ such that for every irreducible component $\partial Y_i$ of $\partial Y$, the reduced inverse image $f^{-1}(\partial Y_i)_{\rm red}$ is a  strict normal crossing divisor on $X$ and satisfies $f(X\setminus |\partial X|) \subseteq Y\setminus |\partial Y|$.   
\end{df}
Note that the composition of morphisms in $\Smlog$ is well defined. Indeed, let $f\colon (X, \partial X)\to (Y, \partial Y)$ and $g\colon (Y, \partial Y) \to (Z, \partial Z)$ be as in Definition \ref{def:Smlog}. It is clear that $g(f(X\setminus |\partial X|)) \subseteq Z\setminus |\partial Z|$, since  $g(f(X\setminus |\partial X|)) \subseteq g(Y\setminus |\partial Y|) \subseteq Z\setminus |\partial Z|$ by assumption. As for the other condition, we need to verify that for every irreducible component $\partial Z_i$ of $\partial Z$, the support of its inverse image $f^{-1}(g^{-1}(\partial Z_i))$ is a strict normal crossing divisor on $X$. Write $g^{-1}(\partial Z_i)_{\rm red} = \cup_{j\in A_i} \partial Y_{ij}$, and $f^{-1}(\partial Y_{ij})_{\rm red} =\cup_{k\in A_{ij}} \partial X_{ijk}$. Thus $f^{-1}(g^{-1}(\partial Z_i))_{\rm red} = \cup_k \cup_j \partial X_{ijk}$ (there might be repetitions in the union, but this is not relevant when considering the support). Since any reduced divisor contained in a    strict normal crossing divisor is itself strict normal crossing, we can conclude.

Suppose now that we are given two pairs $(X,\partial X)$ and $(Y,\partial Y)$ in $\Smlog$. Write  $\partial X_1,\ldots, \partial X_M$ for the components of  $\partial X$, and  $\partial Y_1,\ldots, \partial Y_N$ for the components of $\partial Y$. Their product $(X,\partial X)\times (Y,\partial Y)$ is by definition the pair $(X\times Y, X\times \partial Y + \partial X\times Y)$, where $(X\times \partial Y + \partial X\times Y)$ is by definition the normal crossing divisor formed by $X\times \partial Y_1,\ldots, X\times \partial Y_N, \partial X_1\times Y, \ldots \partial X_M\times Y$. This construction is the categorical product in $\Smlog$. The terminal object in $\Smlog$ is the pair $(\Spec{k}, \emptyset)$.

\begin{df}\label{def:Min_Map_Smlog} Let $f\colon (X, \partial X) \to (Y,\partial Y)$ be a morphism in $\Smlog$. We say that $f$ is \textit{minimal} if $\partial X = f^{-1}(\partial Y)_{\rm red}$.
\end{df}
As a remark, we note that the projection maps $(X\times Y, X\times \partial Y + \partial X\times Y)$ to  $(X, \partial X)$ and $(Y, \partial Y)$ are typically not minimal. In fact, both are minimal if and only if the divisors $\partial X$ and $\partial Y$ are empty. Minimal morphisms play a special role in the construction of certain fiber products, as we will see below. 
\subsubsection{}\label{subsec:adjoint-pair-log-smooth}Let $\omega\colon \Smlog\to \Sm$ be the functor $(X,\partial X)\mapsto X\setminus |\partial X|$. This functor has an obvious left adjoint
\begin{equation}\label{eq:lambda-omega}\lambda\colon \Sm \leftrightarrows \Smlog\colon \omega\end{equation}
that sends a smooth $k$-scheme $X$ to the pair $(X, \emptyset)$. Indeed, any morphism $f\colon (X, \emptyset)\to (Y, \partial Y)$ in $\Smlog$ is given by a morphism of $k$-schemes $f\colon X\to Y$ that has to factor through the open embedding $Y\setminus \partial Y = \omega( (Y, \partial Y)) \to Y$.  Note that the functors $\lambda$ and $\omega$ both commute with products, and that the unit of the adjunction is the identity, expressing $\Sm$ as a retract of $\Smlog$. There is also another functor
\[F\colon \Smlog\to \Sm, \quad P = (X,\partial X)\mapsto F(P) = X,\]
that does not have any obvious left adjoint. If confusion does not arise, given a smooth $k$-scheme $X$, we will write just $X$ in $\Smlog$ for the pair $(X, \emptyset) = \lambda(X)$.
\subsubsection{}\label{ssec:def-cuben} Let $\P^1$ be the projective line over $k$ and let $y$ be the standard rational coordinate on it. For every $n\geq 1$ we have a distinguished object $\boxq^n$ in $\Smlog$, defined as $\boxq^n = ((\P^1)^n, F^n_\infty)$, where $F^n_\infty$ denotes the normal crossing divisor $\sum_{i=1}^n(y_i =\infty)$. There are also maps $\iota_{\varepsilon, i}^n\colon \boxq^n\hookrightarrow \boxq^{n+1}$ for $\varepsilon \in \{0,1\}, i = 1, \ldots n+1$ given by the inclusion with
\[ (\iota_{\varepsilon, i}^{n})^*(y_j) = y_j \text{ for $1\leq j < i$ }, (\iota_{\varepsilon, i}^{n})^*(y_j) = y_{j-1} \text{ for $i< j \leq n+1$ and }    (\iota_{\varepsilon, i}^{n})^*(y_i) = \varepsilon,
\]
as well as projections $p^n_i\colon \boxq^n\to \boxq^{n-1}$ for $i=1,\ldots, n$ induced by $p^n_i\colon (\P^1)^n\to (\P^1)^{n-1}$ that forgets the $i$-th coordinate. We will denote by $\delta_n\colon \boxq^n\to \boxq^n\times \boxq^n$ the diagonal map.
\begin{rmk}To get a feeling of the importance of the additional datum of a ``divisor at infinity'' $\partial X$ on a smooth scheme $X$, consider the pairs $\mathbb{A}^1=(\mathbb{A}^1,\emptyset)$ and $\boxq^1 = (\P^1, \infty)$ and $\P^1 = (\P^1, \emptyset)$. There are canonical maps
    \[ \mathbb{A}^1\to \boxq^1\to \P^1,\]
the first being given by the counit of the adjunction displayed in \eqref{eq:lambda-omega},   and there are no non-costant maps in the opposite directions, so that the three objects are all distinct in $\Smlog$.
    \end{rmk}
\subsection{Modulus pairs \`a la Kahn-Saito-Yamazaki} We recall from \cite{KSY2}\footnote{The paper \cite{KSY2} has been superseded by \cite{KMSY} during the revision of this manuscript. Since the latter is not yet available to the public, we have decided to keep the references to the obsolete version. To the best of our knowledge, the relevant parts will appear unchanged in \cite{KMSY} (but presumably with a different numbering).}  the following constructions.
\begin{df}\label{def:mod-pair-KSY}A modulus pair $M=(\ol{M}, M^\infty)$ consists of a scheme $\ol{M}\in \mathbf{Sch}(k)$ and an effective Cartier divisor  $M^\infty\subset \ol{M}$, possibly empty, such that $\ol{M}$ is locally integral and the open (dense) subset $M^o = \ol{M}\setminus M^\infty$ is smooth and separated over $k$.  A pair $M$ is called proper if $\ol{M}\to \Spec{k}$ is proper. 
\end{df}
\begin{df}\label{def:Adm_map_MSm}Let $M_1, M_2$ be two modulus pairs. Consider a scheme-theoretic morphism \[f\colon M_1^o\to M_2^o\] over $k$. Let $\ol{\Gamma}_f\subset \ol{M}_1\times \ol{M}_2$ be the closure of the graph of $f$ and let $p_1, p_2$ be the two projections
\[p_1\colon \ol{M}_1\times \ol{M}_2 \to \ol{M}_1,\quad p_2\colon \ol{M}_1\times\ol{M}_2 \to \ol{M}_2.\] 
Let $\phi\colon \ol{\Gamma}_f^N\to \ol{M}_1\times \ol{M}_2$ be the composition of the normalization morphism of the closure of the graph with the inclusion. We say that $f$ is admissible for the pair $M_1, M_2$ if
\begin{romanlist}
\item the composition morphism $p_1\circ \phi \colon \ol{\Gamma}^N_f \to \ol{M}_1$  is proper,
\item there is an inequality $\phi^* p_1^*(M^\infty_1) \geq \phi^*p_2^*(M^\infty_2)$ as Weil divisors on $\ol{\Gamma}^N_f$. 
\end{romanlist}
\end{df}
Here we are using the standard partial ordering on the set of Weil divisors on a scheme by setting $D\geq D'$ if and only if $D-D'$ is effective. 
We denote by $\mathbf{\ul{M}Sm}$ the category having as objects modulus pairs and as morphism admissible morphisms between them (the fact that composition of morphisms is well defined is verified in \cite{KSY2}, Section 1.2). With $\mathbf{\ul{M}Sm}^{\rm fin}$, we denote the subcategory of $\mathbf{\ul{M}Sm}$ having the same objects and whose morphisms satisfy the additional condition that  $p_1\circ \phi \colon \ol{\Gamma}^N_f \to \ol{M}_1$ is finite. Finally, we  denote by $\mathbf{MSm}$ the full subcategory of $\mathbf{\ul{M}Sm}$ whose objects are proper modulus pairs.
\begin{rmk} Note the difference between the admissibility condition of a morphism in $\mathbf{\ul{M}Sm}$  and a morphism in $\mathbf{Sm_{log}}$, even in the case a morphism $f\colon M_1\to M_2$ in $\mathbf{\ul{M}Sm}$ is induced by a morphism of smooth schemes $f\colon \ol{M_1}\to \ol{M_2}$ (e.g. if $f$ is a map in $\mathbf{\ul{M}Sm}^{\rm fin}$). In this situation, condition ii) in \ref{def:Adm_map_MSm} reads
\begin{equation}\label{eq:inequality-Cartier-MSm}\nu^* M_1^\infty \geq \nu^* f^* M_2^\infty
\end{equation}
where $\nu\colon \ol{M_1}^N\to \ol{M_1}$ is the normalization morphism. The inequality in \eqref{eq:inequality-Cartier-MSm} is an inequality of effective Weil divisors on a normal variety. Suppose now that both $\ol{M_1}$ and $\ol{M_2}$ are smooth and that the divisors $M_1^\infty$ and $M_2^\infty$ have SNC support. We immediately see that $M_1^\infty \geq f^* M_2^\infty $ implies that $f$ gives rise to a map in $\Smlog$ according to Definition  \ref{def:Smlog}. On the other hand, since the admissibility condition in $\Smlog$ is checked on the reduced pull-back of the divisor, it is not true that a map in $\Smlog$ gives rise to a map of pairs in $\mathbf{\ul{M}Sm}$ (see  Remark \ref{rmk:useful-example} below for a useful example). 
\end{rmk}
\subsubsection{}\label{def:tensor_productMSm} For $M, N \in \mathbf{\ul{M}Sm}$, we define their tensor product $L =M\otimes N$  by $\ol{L} = \ol{M}\times \ol{N}$ and $ L^\infty = M^\infty \times \ol{N} + \ol{M} \times N^\infty$. This gives the categories $\mathbf{\ul{M}Sm}$ and $\mathbf{MSm}$ a symmetric monoidal structure, with unit the modulus pair $(\Spec{k},\emptyset)$.  

\begin{rmk}\label{rmk:useful-example} As noticed in \cite{KSY2}, Warning 1.12, the tensor product $M_1\otimes M_2$ does not have the universal property of products, since, for example, the diagonal morphism $M\xrightarrow{\Delta} M\otimes M$ is not admissible as soon as $M^\infty$ is not empty. Indeed, let $M = (\P^1, \infty)$. Then $M\otimes M = (\P^1\times_k \P^1, \infty\times \P^1 + \P^1\times \infty)$. The diagonal map $\delta\colon \P^1\to \P^1\times \P^1$ is not admissible in $\mathbf{\ul{M}Sm}$, since 
\[\delta^*(\infty\times \P^1 + \P^1\times \infty)  = 2\cdot \infty \nleq \infty.\]
On the other hand, the map $\delta=\delta_1$ is a morphism in $\Smlog$ between $\boxq^1$ and $\boxq^1\times \boxq^1 = \boxq^2$. \end{rmk}
\subsection{Inverting birational maps}\label{subsec:inv-rat-maps} We need to embed the category $\Smlog$ in a larger category with the same objects but where we allow some morphims to be defined only after a proper birational transformation. Recall that we are assuming that $k$ \textit{admits resolution of singularities}, i.e. that the following two conditions hold:
\begin{enumerate}
    \item for any reduced scheme of finite type $X$ over $k$, there exists a proper birational morphism $f\colon \tilde{X}\to X$ such that $\tilde{X}$ is smooth,
    \item for any smooth scheme $X$ over $k$ and a proper surjective morphism $Y\to X$ which has a section over a dense open subset $U$ of $X$, there exists a sequence of blow-ups  $X_n\to X_{n-1}\to\ldots X_0 = X$, with smooth centers lying over $X\setminus U$, and a morphism $X_n\to Y$ such that the composition $X_n\to Y\to X$ is the structure morphism $X_n\to X$. 
    \end{enumerate}

\begin{df}Let $P = (X, \partial X)\in \Smlog$. We denote by $\mathfrak{B}_P$ the category of \textit{admissible blow-ups} of $P$. An object of $\mathfrak{B}_P$ is a morphism $\pi\colon P'=(X', \partial X')\to P$ with $P' \in \Smlog$ induced by a projective birational map $\pi\colon X'\to X$ that restricts to an isomorphism on the open complements $\omega(X')\cong \omega(X)$ and such that $|\partial X'| = |\pi^{-1}(\partial X)|$. Morphisms in $\mathfrak{B}_P$ are the minimal morphisms in $\Smlog$ over $P$. If  $\mathcal{S}_b$ denotes the class of admissible blow-ups of pairs, $\mathfrak{B}_P$ is the full subcategory of the comma category $\Smlog/P$ given by the objects $P'\xrightarrow{s} P$ with $s \in \mathcal{S}_b$.
\end{df}
In other words, an object in $\mathfrak{B}_P$ is a blow-up with center in a closed subscheme $Z\subset \partial X\subset X$.

\begin{prop}\label{prop:calc-fractions}The class $\mathcal{S}_b$ enjoys a calculus of right fractions. In particular, for every $P,Q\in \Smlog$, the natural map
    \[  \varcolim_{P'\in \mathfrak{B}_P} \hom_{ \Smlog }(P', Q)\to \hom_{\Smlog[\mathcal{S}_b^{-1}]} (P,Q)\]
    is an isomorphism. Moreover, since for any $P\in \Smlog$ the category $\mathfrak{B}_P$ contains a small cofinal subcategory,  the ${\rm Hom}$ sets of $\Smlog[\mathcal{S}_b^{-1}]$ are small.
    \begin{proof}We recall from \cite{GZFractions} the conditions that a class of morphisms $\Sigma$ has to satisfy in order to enjoy calculus of right fractions (this is dual to \cite[I.2.2]{GZFractions}). These conditions are: a) the identities of $\Smlog$ are in $\Sigma$; b)   if $u\colon X\to Y$ and $v\colon Y\to Z$ are in $\Sigma$, then their composition $v\circ u$ is also in $\Sigma$; c) for each diagram $X'\xrightarrow{s} X\xleftarrow{u} Y $ where $s\in \Sigma$, there exists a commutative square 
        \begin{equation}\label{eq:orecondition-calc-frac}\xymatrix{ Y' \ar[r]^{u'} \ar[d]^{s'}  & X'\ar[d]^s \\ 
        Y\ar[r]^u & X
          }\end{equation}
with $s'\in \Sigma$; d) if $f,g\colon X\rightrightarrows Y$ are morphisms of $\Smlog$ and if $s\colon Y\to Y'$ is a morphism of $\Sigma$ such that $s\circ f =  s\circ g$, then there exists a morphism $t\colon X'\to X$ such that $f\circ t = g\circ t$. For $P\in \Smlog$, write $\Sigma \downarrow P$ for the full subcategory of the comma category $\Smlog/P$ given by the objects $P'\xrightarrow{s} P$  with $s\in \Sigma$. If a class of arrows $\Sigma$ enjoys calculus of right fractions, there is an isomorphism $\colim_{P'\in \Sigma \downarrow P} \hom_{\Smlog}(P', Q)\xrightarrow{\cong} \hom_{\Smlog[\Sigma^{-1}]} (P, Q)$, natural in $P$ and $Q$, by the dual of \cite[Proposition I.2.4]{GZFractions}. Note that by resolution of singulartities, for each $(X, \partial X)\in \Smlog$, the category $\mathfrak{B}_{(X,\partial X)}$ contains the the cofinal subcategory consisting of maps $(X_n, \partial X_n) \to (X, \partial X)$ where $X_n\to X$ is a composition of blow-ups with smooth centers contained in $X\setminus |\partial X|$ (the cofinality follows directly from condition (2) in the above definition of resolution).
 Let's prove that  $\mathcal{S}_b$ enjoys calculus of right fractions. Conditions a) and b) are obvious. We first show that condition c) holds. Given pairs $(X, \partial X)$, $(X', \partial X')$ and $(Y,\partial Y)$ with maps $f\colon (Y, \partial Y) \to (X, \partial X)$ and $\pi \colon (X', \partial X')\to (X, \partial X)$ with $\pi \in \mathfrak{B}_{(X, \partial X)}$, we define the pair $(Y', \partial Y')$ as follows. Set $\tilde{Y} = Y\times_X X'$ and write $p_Y$ for the projection $\tilde{Y}\to Y$. By assumption, $\tilde{Y}\setminus| p_Y^{-1}(f^{-1}(\partial X))|$ is isomorphic to $Y\setminus |f^{-1} \partial X|$. By resolution of singularities, we can find a projective birational map 
 \[\pi'\colon Y'\to \tilde{Y} \to Y\] 
 that is obtained as sequence of blow-ups with smooth centers lying over $\partial Y$ and such that $\partial Y' := (\pi')^{-1}(\partial Y)_{\rm red}$ is a normal crossing divisor on $Y'$. Then $(Y', \partial Y') \to (Y, \partial Y)$ is a morphism in $\mathfrak{B}_{(Y,\partial Y)}$. Note that the induced morphism $f'\colon (Y',\partial Y')\to (X', \partial X')$ is  admissible and gives a commutative square like \eqref{eq:orecondition-calc-frac} as required. The commutativity is clear, so we only need to verify the admissibility conditions of Definition \ref{def:Smlog}. Since $f$ is admissible, we have $f(Y \setminus |\partial Y|) \subseteq X \setminus |\partial X| = X' \setminus |\partial X'|$ (where the latter equality follows from the fact that $\pi$ is in $\mathfrak{B}_{(X, \partial X)}$). Since $\pi'$ belongs to $\mathfrak{B}_{(Y, \partial Y)}$, it follows that $Y'\setminus |\partial Y'| = Y\setminus |\partial Y|$, and the induced morphism $f'$ restricts to $f$ on the open complements. In particular, we have that $f'(Y' \setminus |\partial Y'|) \subseteq X' \setminus |\partial X'|$. Let now $\partial X'_i$ be a component of $\partial X'$. Since $\partial X'= \pi^{-1}(\partial X)_{\rm red}$, there is a component $\partial X_j$ of $\partial X$ such that $\partial X'_i \subseteq \pi^{-1}(\partial X_j)_{\rm red}$. Then $(f')^{-1}(\partial X_i')_{\rm red} \subseteq (\pi')^{-1}(f^{-1}(\partial X_j))_{\rm red}\subseteq \partial Y'$, which is a strict normal crossing divisor on $Y'$ by construction. So $f'$ is admissible. 

Finally, we turn to condition d). Suppose we are given $f,g\colon (X, \partial X) \to (Y, \partial Y)$ and a morphism $s\colon (Y, \partial Y) \to (Y', \partial Y)$ in $\mathcal{S}_b$ such that $s\circ f = s\circ g$. By assumption, $s$ is given by a projective birational morphism $Y\to Y'$ which restricts to the identity on the open complements $\omega((Y, \partial Y)) = \omega((Y', \partial Y'))$. In particular, the underlying morphisms $\omega(f), \omega(g)\colon \omega((X, \partial X)) = X\setminus |\partial X| \to Y\setminus |\partial Y|$ agree. The required morphism $t$ is then given by the counit map $\lambda(X\setminus \partial X) \to (X,\partial X)$, noting that the compositions $f\circ t$ and $g\circ t$ have to factor through $\lambda(Y\setminus |\partial Y|) \to (Y, \partial Y)$. 
\end{proof}
    \end{prop}
    
Write $\BSmlog$ for the localized category $\Smlog[\mathcal{S}_b^{-1}]$. We denote by 
\begin{equation}\label{eq:loc-functorSmlog}v\colon \Smlog \to \BSmlog\end{equation} the localization functor: it is  clearly faithful (note that admissible blow-ups are in particular epimorphisms of schemes), and by \cite[I.3.6]{GZFractions} commutes with finite direct and inverse limits that exist in $\Smlog$. Note here that finite products exist in $\Smlog$ and there is a terminal object $(\Spec{k}, \emptyset)$, but    arbitrary fiber products are not representable in $\Smlog$: indeed, let $X, X'$ be smooth hypersurfaces in $Z=\P^n_k$ such that their intersection is not smooth. We can consider $X, X'$ and $Z$ as objects of $\Smlog$ with empty boundary (via $\lambda$). Since $\omega$ is a right adjoint, it preserves finite limits, and therefore if the product $\lambda(X)\times_{\lambda(Z)} \lambda(X')$ existed in $\Smlog$, so would the product $\omega(\lambda(X)\times_{\lambda(Z)} \lambda(X')) = X\times_Z X'$ in $\Sm$. Since $X\times_Z X'$ is not smooth, this is not the case. This implies by \cite[Tag 04AS]{stacks-project}, that $\Smlog$ does not have all small limits. 
	
We will construct in \ref{cor:quarrable_smooth_morphisms} fiber products in $\Smlog$ where one of the maps is a smooth morphism that is minimal in the sense of Definition \ref{def:Min_Map_Smlog}). 

\subsubsection{} If an object in $\Smlog$ has empty boundary divisor, our definition does not allow more morphisms to appear in $\BSmlog$. More precisely, we have
\[\hom_{\Smlog}((X, \emptyset), (Y, \partial Y)) \xrightarrow{\cong}\hom_{\BSmlog}((X, \emptyset), (Y, \partial Y)) \]
for every $(Y, \partial Y)$ in $\Smlog$. As a consequence of this fact, we can extend the adjunction of \ref{subsec:adjoint-pair-log-smooth} to the localized category $\BSmlog$. Indeed, we first notice that  every admissible blow-up $\mathfrak{B}_P$ for a given pair $P = (X, \partial X)$ does not change the open complement $X\setminus \partial X$ (up to isomorphism). The universal property of the localization allows then to define a functor
\[\omega\colon \BSmlog \to \Sm, \quad (X,\partial X)\mapsto X\setminus \partial X,\]
that restricts to the functor $\omega$ of \ref{subsec:adjoint-pair-log-smooth} when composed with the localization functor $v$. The above observation shows then that we have a bijection
\[\hom_{\Sm}(X, (Y\setminus \partial Y)) \xrightarrow{\cong}\hom_{\BSmlog}((X,\emptyset), (Y, \partial Y))\]
so that the composite functor $\lambda\colon \Sm \to \Smlog \xrightarrow{v}\BSmlog$ is left adjoint to $\omega$.  

A similar situation shows up in the case $\dim X = 1$. In this case, morphisms from $(X,\partial X)$ to any other object of $\Smlog$ do not change if we pass to $\BSmlog$, since every morphism in $\mathfrak{B}_{(X, \partial X)}$ is an isomorphism already in $\Smlog$.
\begin{rmk}Let $\PSh({\Smlog})$ (resp. $\PSh(\BSmlog)$) be the category of presheaves of sets on $\Smlog$ (resp. on $\BSmlog$). Then the functor $v$ induces a string of adjoint functors $(v_!, v^*, v_*)$ (where each functor is the left adjoint to the the following one) between the categories of presheaves: $v^*$ is induced by composition with $v$, while $v_!$ (resp.~$v_*$) is the left (resp.~right) Kan extension of $v$. Since $v$ is a localization, then $v_!$ is also a localization or, equivalently (by \cite[Proposition I.1.3]{GZFractions}), $v^*$ is fully faithful. The functor $v^*$ identifies $\PSh(\BSmlog)$ with the subcategory of $\PSh(\Smlog)$ of presheaves that invert the morphisms in $\cS_b$.
\end{rmk}
\subsection{Modulus data} We now introduce the category of \textit{modulus data} over $k$, that will be the basic object for our constructions.
\begin{df}\label{def:modulus-data}A \textit{modulus datum} $M$ consists of a triple $M = (\ol{M}; \partial M, D_M )$, where $\ol{M}\in \Sm$ is a smooth $k$-scheme, $\partial M$ is a strict normal crossing divisor on $M$ (possibly empty), $D_M$ is an effective Cartier divisor on $M$ (again, the case $D_M = \emptyset$ is allowed), and \textit{the total divisor $|D_M|_{\rm red} + \partial M$ is a strict normal crossing divisor on $\ol{M}$}.
    
    Let $M_1, M_2$ be two modulus data. A morphism $f\colon M_1\to M_2$ is called \textit{admissible} if it is a morphism of $k$-schemes $f\colon \ol{M_1}\to \ol{M_2}$ that satisfies the following conditions.
    \begin{romanlist}
        \item The map $f$ is a morphism in $\Smlog$ between $(\ol{M_1}, \partial{M_1})$ and $(\ol{M_2}, \partial{M_2})$, i.e. for every irreducible component $\partial M_{2,l}$ of $\partial M_2$ we have $|f^{*}(\partial M_{2,l})|\subseteq  |\partial M_{1}|$. 
        \item The divisor $f^*(D_{M_2})$ is defined, and satisfies $f^*(D_{M_2}) \geq D_{M_1}$ as Weil divisors on $\ol{M_1}$.
        \end{romanlist}
    We denote by $\MSmlog$ the category having objects modulus data and morphisms admissible morphisms. If one inverts the inequality in condition ii)  above, we obtain a ``dual'' category $\mathbf{\ul{M}}\Smlog$. We will refer to condition ii) as the \textit{modulus condition} on morphism. If equality holds in ii), we will say that the morphism $f$ is \textit{minimal with respect to the modulus condition}. If $f$ is minimal also with respect to the boundary divisors $\partial M_1$ and $\partial M_2$ in the sense of Definition \ref{def:Min_Map_Smlog}, we will simply say that $f$ is a minimal morphism of modulus data. Finally, given a modulus datum $M=(\ol{M}; \partial M, D_M)$, we will say that $\partial M$ is the \textit{boundary divisor} of the datum $M$ and that $D_M$ is its \textit{modulus divisor}.
    \end{df}
    \begin{rmk}The condition that  $|D_M|_{\rm red} + \partial M$ forms a strict normal crossing divisor on $\ol{M}$ is \textit{not} necessary for the construction of our motivic homotopy category $\MH$. This condition is however useful for the construction of a $K$-theory space associated to a modulus datum $M$, as discussed in \cite{BThesis}.
        \end{rmk}
\begin{df}\label{def:mod-pairs}The category of \textit{modulus pairs} $\mathbf{\ol{M}Sm}(k)$ is the category having objects pairs $(M, D_M)$ where $M$ is a smooth $k$-scheme and $D_M$ is an effective Cartier divisor on it such that the open complement $M^o= M\setminus |D_M|$ is dense.   A morphism of pairs $M_1\xrightarrow{f} M_2$ is a morphism of $k$-schemes such that $f^*(D_{M_2}) \geq D_{M_1}$ as Weil (or Cartier) divisors on $M_1$.
    \end{df}
    \begin{rmk}We are here using the opposite inequality of the definition given in \cite{KSY2}. A possible way for unifying the two notions would be to allow \textit{non effective} modulus pairs $M=(\ol{M}, D_M)$ in $\mathbf{\ul{M}Sm}$ where $D_M$ is a Cartier divisor on $\ol{M}$, not necessarily effective. Then one can fix the direction of the inequality of the definition of modulus pairs as done in \cite{KSY2}, Definition 1.1, and embed our category $\mathbf{\ol{M}Sm}(k)$  in $\mathbf{\ul{M}Cor}$ by sending $(\ol{M}, D_M)$ to $(\ol{M}, -(D_M))$.
        \end{rmk}
\subsubsection{}\label{subsubsec:vfunctor}There is a fully faithful functor $u\colon \Smlog\to \MSmlog$, that sends a pair $P = (X, \partial X)$ to the modulus datum  $u(P)=(X; \partial X, \emptyset)$, as well as a ``forgetful'' functor $F\colon \MSmlog\to \mathbf{\ol{M}Sm}$ that sends a modulus datum $(\ol{M};\partial M, D_M)$ to the modulus pair $(\ol{M}, D_M)$. They fit together in a commutative square of categories
\[\xymatrix{ \MSmlog \ar[r]^{F}& \mathbf{\ol{M}Sm}(k),\\
\Smlog\ar[u]^u \ar[r]^{F} & \Sm\ar[u]^u
}\quad
\xymatrix{(\ol{M}; \partial M, D_M) \ar@{|->}[r] & (\ol{M}, D_M)\\
(X,\partial X) \ar@{|->}[r] & X
}
\]
\subsection{Fiber products} As for products, fiber products do not exist in general in the categories $\Smlog$ and $\MSmlog$. We have, however, the following useful proposition.
\begin{prop}\label{prop:quarrable_smooth_morphisms} Let $f\colon M = (\ol{M};\partial M, D_M)\to N = (\ol{N}; \partial N, D_N)$ be a minimal morphism in $\MSmlog$ such that the underlying morphism of schemes $f\colon \ol{M}\to\ol{N}$ is smooth. Then for every $g\colon L = (\ol{L};\partial L, D_L)\to N$, the fiber product $L\times_N M$ exists in $\MSmlog$.
    \begin{proof}Since $f$ is smooth, the fiber product $\ol{M}' = \ol{L}\times_{\ol{N}} \ol{M}$ is also smooth over $k$. Write $f'\colon \ol{M}'\to L$ for the base-change map, and define $\partial{M}'$ to be the divisor $(f')^{-1}(\partial L)_{\rm red}$. Each component of $\partial M'$ is the inverse image of a (smooth) component of $\partial L$ along a smooth map, so it is smooth over $k$. Similarly, each face $\partial M'_I = \cap_{i\in I}\partial M'_i$ is smooth over $k$ and of pure codimension $|I|$ on $\ol{M}'$. Thus, $\partial M'$ defined in this way is a strict normal crossing divisor on $\ol{M}'$. By construction, it's clear that the morphisms $f'$ and $g'\colon \ol{M}'\to \ol{M}$ are admissible morphisms in $\Smlog$, and that $f'$ is minimal. As for the modulus condition, set $D_{M'}$ to be the divisor $(f')^*(D_L)$. Then we have
        \[ (g')^*(D_M) = (g')^* f^* D_N = (f')^* g^* D_N \geq (f')^*  D_{L} = D_{M'},\]
        where the first equality follows from the minimality requirement on $f$, so that the maps $g'$ and $f'$ are both satisfying the modulus condition, and therefore are admissible morphisms in $\MSmlog$. We are  left to show that the universal property of the fiber product is satisfied by the modulus datum $M'$, but this is straightforward.
        \end{proof}
    \end{prop}
We deduce from the case  $D_M= D_L = \emptyset$ the analogous statement for $\Smlog$.
\begin{cor}\label{cor:quarrable_smooth_morphisms} Let $f\colon (X, \partial X) \to (Y,\partial Y)$ be a minimal morphism in $\Smlog$ such that $f\colon X\to Y$ is smooth. Then for every map $g\colon (Z, \partial Z) \to (Y,\partial Y)$, the fiber product $(X, \partial X)\times_{(Y,\partial Y)} (Z, \partial Z)$ is representable in $\Smlog$. 
    \end{cor}
Of course, there is nothing to say in case $\partial X$ and $ \partial Y$ are also empty.
\subsection{Monoidal structure on $\MSmlog$}\label{subsec:monoidal-structure-MSmlog} We extend the product in $\Smlog$ to a symmetric monoidal structure on $\MSmlog$.
\begin{df}Let $M,N\in \MSmlog$ be modulus data. We define the modulus datum $L = M\tensor N$ by
    \[L = (\ol{M}\times \ol{N}; \partial L = \partial M \times\ol{N} +  \ol{M}\times \partial M, D_L = D_M\times \ol{N} + \ol{M}\times D_N).\]
    The category $\MSmlog$ equipped with the tensor product $\tensor$ is a symmetric monoidal category, with unit object $ \mathds{1} = (\Spec{k}, \emptyset, \emptyset)$. In a similar fashion, we define a symmetric monoidal product $\tensor$ on the category of modulus pairs  $\mathbf{\ol{M}Sm}(k)$. 
    \end{df}
\subsubsection{}\label{subsubsec:tensor-cart-empty-modulus}When $M=(\ol{M}; \partial M, \emptyset)$ and $N = (\ol{N};\partial N, \emptyset)$, then we can check that $M\tensor N = M\times N = (\ol{M}\times \ol{N}; \partial L, \emptyset)$ is the categorical product of $M$ and $N$. In particular, the functor $u\colon \Smlog\to \MSmlog$ is strict monoidal (when one considers on $\Smlog$ the monoidal structure given by the cartesian product). The forgetful functor $F\colon \MSmlog\to  \mathbf{\ol{M}Sm}(k)$ is also strict monoidal.
\begin{rmk}For arbitrary objects $M,N$ in $\MSmlog$, our choice of admissibility condition for morphisms prevents the existence of projection maps $M\tensor N\to M$ or $M\tensor N\to N$. However, if $M= (\ol{M}; \partial M, \emptyset)$ and $N= (\ol{N};\partial N, D_M)$, then the map $M\tensor N \to N$ induced by the projection $\ol{M}\times \ol{N}\to \ol{N}$ is clearly admissible.  
    \end{rmk}
\subsection{A digression on interval objects in monoidal categories}\label{subsec:Inverval-objects-monoidalcat} Let $(\mathcal{M}, \tensor, \mathds{1})$ be a symmetric monoidal category with unit object $\mathds{1}$. 
\begin{df}An object $I$ in $\mathcal{M}$ is called a \textit{weak interval} in $\mathcal{M}$ if there exist a map $p_I\colon I \to \mathds{1}$ (the ``projection'') and monomorphisms $\iota^I_{\varepsilon}\colon \mathds{1}\to I$ for $\varepsilon = 0,1$ (the ``inclusions at $0$ and $1$'') that satisfy
    \[p_I\circ \iota^I_0 = p_I\circ \iota^I_1 =\id_{\mathds{1}}.\]
\end{df}
Let $\mathcal{M}$ be a symmetric monoidal category equipped with a weak interval $I$.
An $\Itensor-$homotopy between two maps $f, g\colon X\rightrightarrows Y$ is the datum of a morphism $H\colon X\tensor I \to Y$ in $\mathcal{M}$ such that $f = H\circ (\id_X\tensor \iota_0^I)$  and $g = H\circ (\id_X\tensor \iota_1^I)$.
\begin{df}\label{def:int-obj-mon-cat}An object $I$ in $\mathcal{M}$ is called an \textit{interval} in $\mathcal{M}$ if it is a weak interval $(I, \iota^I_0, \iota^I_1, p_I)$ that is additionally equipped with a multiplication map
    \[\mu\colon I\tensor I\to I,\]
    verifying the identities $\mu\circ (\id_I\tensor \iota^I_0) = \iota_0^I \circ p_I$ and $\mu\circ (\id_I\tensor \iota^I_1) = \id_{I}$.
    \end{df}
The notion of interval object presented here is more general then the definition of Voevodsky in \cite{VSelecta}.

\subsubsection{}Any weak interval object $I$ in $\mathcal{M}$ determines a co-cubical object (see \cite{LevineSmooth} for the notation) $I^\bullet\colon \mathbf{Cube} \to \mathcal{M}$ by setting
\[I^n = I^{\tensor n},\quad p^n_{i, I} = \id_I^{\tensor (i-1)}\tensor p_I\tensor \id_I^{\tensor (n-1)}, \quad \delta^n_{i,\varepsilon} =  \id_I^{\tensor (i-1)}\tensor \iota^I_\varepsilon \tensor \id_I^{\tensor (n-1)}  \]
for $\varepsilon\in\{0,1\}$, $i = 1,\ldots, n$. If $I$ is moreover an interval object (so that is equipped with a multiplication map), the same formulas work to give an extended co-cubical object, $I^\bullet\colon \mathbf{ECube} \to \mathcal{M}$, where $\mathbf{ECube}$ is the extended cubical category introduced e.g.~in \cite{LevineSmooth}.

Conversely, given an extended cocubical object $C$ in $\mathcal{M}$, where the monoidal structure on $\mathbf{ECube}$ is given by cartesian product, one can easily check that $I=C([1])$ is an interval object in $\mathcal{M}$.
 
\subsubsection{} In arbitrary monoidal categories there are no diagonal morphisms,  so that given a weak interval object one can~---~a priori~---~only develop a cubical theory and not a  simplicial theory. Fortunately,  we will consider for our applications an interval object $(I, \iota_1^I, \iota_0^I, p_I, \mu)$ that is equipped with an extra map $\delta_I\colon I\to I^{\tensor 2}$ such that the compositions $(\id_I\tensor p_I)\circ \delta$ and $( p_I\tensor \id_I)\circ \delta$ are the identity on $I$.

Following \cite[2.2, p.118-119]{VSelecta}, we can then construct a universal cosimplicial object  in $\mathcal{M}$ as follows. 
\begin{df}\label{def:DeltaI} Set $\Delta_I^n = I^{\tensor n}$ for every $n$. For $i=0,\ldots, n$, let  \[d^{i}\colon [n-1] = \{0,\ldots, n-1\}\to [n] = \{0,\ldots, n\}\quad ({\text{resp. }} s^i\colon [n+1]\to [n]) \] be the standard $i-$th face (resp. $i-$th degeneracy) in the  simplicial category $\Delta$. Define
 \begin{equation}\label{eq:formulaDeltaId}\Delta_I(d^i) = 
 \begin{cases} 
     \iota_0^I\tensor \id_I^{\tensor (n-1)} & \mbox{if } i=0, \\
      \id_I^{\tensor (n-1)} \tensor \iota_1^I & \mbox{if } i=n, \\ 
      \id^{\tensor (i-1)}_I\tensor\delta_I\tensor \id_I^{\tensor(n-i-1)} & \mbox{if } 1\leq i\leq n-1
\end{cases}\end{equation}
and 
\begin{equation}\label{eq:formulaDeltaIs}\Delta_I(s^i)   = \id_I^{\tensor i} \tensor p_I\tensor \id_I^{\tensor (n-i)}.\end{equation} 
\end{df}
    
It is easy to check that this data define a cosimplicial object in $\mathcal{M}$, that we will denote by $\Delta_I^\bullet$. 
\begin{rmk} The formulas in \eqref{eq:formulaDeltaId} and \eqref{eq:formulaDeltaIs} are not explicit in \cite{VSelecta}. Voevodsky's construction holds in a  $\tensor$-category that is a site with products, the tensor structure being given by cartesian products of objects, and this fact is used in the formulation of \textit{loc.cit.}. It is not hard (but a bit tedious) to deduce from Voevodsky's definition our formulas.
\end{rmk}
\subsection{A distinguished interval in $\Smlog$}\label{sec:def-I-object-sec1} We specialize the result of the previous subsection to our case of interest. 

Consider the object $\boxq = \boxq^1 = (\P^1, \infty)$ in $\Smlog$. We have two distinguished admissible morphisms in $\Smlog$
\[\iota_0^{\boxq}, \iota_1^{\boxq}\colon \Spec{k} = (\Spec{k}, \emptyset) \rightrightarrows \boxq\]  
induced by the inclusions of the points $0=[0:1]$ and $1=[1:1]$ in $\P^1$. There is also a projection
\[p_{\boxq}\colon \boxq\to (\Spec{k}, \emptyset) =: \mathds{1}, \]
indued by the structure map $\P^1\to \Spec{k}$ and satisfying the obvious property that \[p_{\boxq}\circ \iota_0^{\boxq} = p_{\boxq}\circ \iota_1^{\boxq }= \id_{\mathds{1}},\]
making $(\boxq, \iota_0^{\boxq },  \iota_1^{\boxq }, p_{\boxq})$ a weak interval object in $\Smlog$.
 
Let $\BSmlog$  as in \ref{subsec:inv-rat-maps} be the localization of the category $\Smlog$ to admissible blow-ups, and consider $\boxq$ as object there. Let $F_\infty$ denote the divisor $(\infty\times \P^1 + \P^1\times \infty)$ on $\P^1\times \P^1$ (we drop the superscript ``2'' from the notation used in \ref{ssec:def-cuben}). There is an extra multiplication map \[\mu\colon (\P^1\times \P^1, F_\infty)\to \boxq\] in $\BSmlog$ induced by the following diagram:
\begin{equation}\label{eq:mult-P1}\xymatrix{ \P^1\times \P^1 \ar@{-->}[r] & \P^1. \\
{\mathbf{Bl}}_{0\times \infty 
+ \infty \times 0}(\P^1\times \P^1) \ar[u]^{\pi} \ar[ru]^{\tilde{\mu}}
}\end{equation}
Here, $\pi\colon B= \mathbf{Bl}_{0\times \infty, \infty \times 0}(\P^1\times \P^1) \to (\P^1\times \P^1)$ is the blow-up along the closed subscheme $(0\times \infty \cup \infty \times 0) \subset F_\infty$. It is easy to see that $B$ is smooth over $k$ (since it is the blow-up along a regularly embedded subscheme), and that it agrees with the closure in $\P^1\times \P^1\times \P^1$ of the graph of the rational map $\mu\colon \P^1\times \P^1\dashrightarrow \P^1$ given by the multiplication map $\mu\colon \A^1\times \A^1\to \A^1,$ $(x,y)\mapsto xy$ (and we will constantly use this identification in what follows). The map $\tilde{\mu}$ in \eqref{eq:mult-P1} is then identified with the composition $B\hookrightarrow (\P^1)^3\xrightarrow{p_{3}} \P^1$, where $p_3$ is the projection to the third factor.

On $B$ we have two distinguished divisors, that we denote by $\tilde{F}_\infty$ and $E$ respectively. The divisor $\tilde{F}_\infty = \tilde{F}_{\infty,1} + \tilde{F}_{\infty,2}  = \P^1\times \infty\times \infty + \infty\times \P^1\times \infty$ is the strict transform of the boundary divisor $F_\infty$ on $\P^1\times \P^1$ along the map $\pi$. The divisor $E=E_{1} + E_{2} = 0\times \infty \times \P^1 + \infty \times 0\times \P^1$ is the exceptional divisor of the blow-up. Together, $( \tilde{F}_{\infty,1}, \tilde{F}_{\infty,2}, E_{1} , E_{2})$ form a strict normal crossing divisor on $B$, that we simply denote by $\partial B$. Then the pair $(B, \partial B)$ is a well-defined object in $\Smlog$ and the map $\tilde{\mu}$ is an admissible morphism in $\Smlog$. Note that we clearly have $\pi^{-1}(F_\infty) = \partial B$, so that $\pi$ is a minimal morphism in $\Smlog$, that is therefore an admissible blow-up for $(\P^1\times \P^1, F_\infty)$ (and hence becomes invertible in $\BSmlog$). To conclude, we have constructed a well-defined morphism 
\[ \mu \colon (\P^1\times \P^1, F_\infty)\to \boxq, \quad \text{ in } \BSmlog\]
as required.

To show that $\boxq$ with this multiplication morphism defines an interval object in $\BSmlog$, we still need to check that the axioms of Definition \ref{def:int-obj-mon-cat} hold. First, note that there is a natural monoidal structure on $\BSmlog$ that makes the localization functor $v$ strict monoidal. Namely, given pairs $(X, \partial X)$ and $(Y, \partial Y)$, define
\begin{equation}\label{eq:mon-structure-BSmlog}(X, \partial X)\tensor(Y, \partial Y) = (X\times Y, \partial X \times Y + X\times \partial Y).
\end{equation}
 This assignment coincides with the cartesian product in $\Smlog$. In fact, more is true: 
 \begin{lem}The product $(X, \partial X)\tensor(Y, \partial Y) = (X\times Y, \partial X \times Y + X\times \partial Y)$ is the categorical product in $\BSmlog$.
	 \begin{proof}We can assume that $X$ and $Y$ are connected. Let $(Z, \partial Z)\in \Smlog$, and suppose we are given two maps $(Z, \partial Z) \to (X, \partial X)$ and $(Z, \partial Z) \to (Y, \partial Y)$ in $\BSmlog$. By Proposition \ref{prop:calc-fractions}, such maps can be described as zig-zags 
		 \[  (Z, \partial Z) \xleftarrow{\pi_f}(Z_f, \partial Z_f) \xrightarrow{f}(X, \partial X),\]
		 \[  (Z, \partial Z) \xleftarrow{\pi_g}(Z_g, \partial Z_g) \xrightarrow{g}(Y, \partial Y),\]
		 where $f$ and $g$ are morphisms in $\Smlog$, and the underlying morphisms to $\pi_f$ and $\pi_g$ are compositions of blow-ups with smooth centers lying over $\partial Z$. In particular, $Z_f\setminus | \partial Z_f| = Z_g \setminus |\partial Z_g| = Z\setminus |\partial Z|$, and $\pi_f$ and $\pi_g$ are minimal morphisms in $\Smlog$. Let $\tilde{Z}$ be the fiber product $Z_f\times_Z Z_g$, and let $\pi\colon Z_{f,g}\to Z_f\times_Z Z_g$ be a resolution of singularities such that the compositions $Z_{f,g}\to Z_f$ and $Z_{f,g}\to Z_g$ are given by compositions of blow-ups with smooth centers, and such that the support of the total transforms of $\partial Z_f$ and $\partial Z_g$ are strict normal crossing divisors on $Z_{f,g}$. Let $\partial Z_{f,g}$ be  the reduced inverse image of $\partial Z$ under the composition $Z_{f,g}\to Z_f \to Z$ (or, equivalently, $Z_{f,g}\to Z_g \to Z$), and consider $(Z_{f,g}, \partial Z_{f,g})$ as an object of $\Smlog$. Note that the projection $(Z_{f,g}, \partial Z_{f,g})\to (Z, \partial Z)$ belongs to $\mathfrak{B}_{(Z, \partial Z)}$. Write $f$ and $g$ for the induced morphisms $Z_{f,g}\to X$ and $Z_{f,g}\to Y$ respectively. This induces a map $h\colon Z_{f,g}\to X\times Y$, and we need to verify that it is admissible in the sense of Definition \ref{def:Smlog}. Every irreducible component of $(\partial X\times Y + X\times \partial Y)$ is clearly of the form $\partial X_i\times Y$ or $X\times \partial Y_j$, thus their pullback to $Z_{f,g}$ via $h$ is obtained as a pullback of $\partial X_i$ (resp.~$\partial Y_j$) via $f$ (resp.~via $g$). Since $f$ and $g$ are admissible to begin with, and since the projections $Z_{f,g}\to Z_f$ and $Z_{f,g}\to Z_g$ are also admissible, we conclude that the first admissibility condition is satisfied. Next, we need to check that $h( Z_{f,g} \setminus |\partial Z_{f,g}|) \subseteq X\times Y \setminus | \partial X\times Y + X\times \partial Y |$. Note that $Z_{f,g} \setminus |\partial Z_{f,g}| = Z\setminus |\partial Z|$, and by assumption $f(Z\setminus |\partial Z|) \subseteq X\setminus |\partial X|$ and $g(Z\setminus |\partial Z|) \subseteq Y\setminus |\partial Y|$. Now, if we let $p_X\colon X\times Y\to X$ and $p_Y\colon X\times Y\to Y$ be the two projections, we have $p_X^{-1} (X\setminus |\partial X|) \supseteq h( Z_{f,g}\setminus |\partial Z_{f,g}|)$ and $p_Y^{-1} (Y\setminus |\partial Y|) \supseteq h( Z_{f,g}\setminus |\partial Z_{f,g}|)$, so that $p_X^{-1} (X\setminus |\partial X|)\cap p_Y^{-1} (Y\setminus |\partial Y|) \supseteq h( Z_{f,g}\setminus |\partial Z_{f,g}|)$, that is precisely what we needed to show. 
		 
		 The above construction determines a unique map in $\BSmlog$ from $(Z, \partial Z)$ to $(X\times Y, \partial X \times Y + X\times \partial Y)$ such that the composition with the projections to $(X, \partial X)$ and to $(Y, \partial Y)$ agree with the original maps. This shows that $(X, \partial X)\tensor(Y, \partial Y)$ satisfies the universal property for being the categorical product. 
	\end{proof}
	 \end{lem}
Consider now the inclusions at $0$ and $1$ of $\boxq$ in $\Smlog$. First, we note that the morphisms 
\[\id_{\boxq}\times \iota_1^{\boxq} \colon \boxq\times \mathds{1} \cong \boxq \to \P^1\times \P^1, \quad  \iota_1^{\boxq}\times \id_{\boxq} \colon  \mathds{1}\times  \boxq \cong \boxq \to \P^1\times \P^1,  \]
automatically factors through   $\pi$, since their image is disjoint from the center of the blow-up. Explicitly, we have the morphism $ \iota_1^{\boxq}\times \id_{\boxq}\colon \boxq\to B$ given by the diagonal embedding $1\times \Delta_{\P^1}\hookrightarrow B$ induced by $x\mapsto (1, x, \mu(1,x)=x)$, and   $\id_{\boxq}\times \iota_1^{\boxq} $ given by the ``twisted'' diagonal embedding $x\mapsto (x, 1, \mu(x,1)=x)$. These maps are clearly admissible in $\Smlog$. Since the morphism $\tilde{\mu}$ is induced by the third projection, we immediately see that we have identities
\[\tilde{\mu} \circ (\id_{\boxq}\times \iota_1^{\boxq}) = \tilde{\mu}\circ(\iota_1^{\boxq}\times \id_{\boxq}) = \id_{\boxq} \quad \text{ in } \Smlog,\]
that descend to the corresponding identities in $\BSmlog$ once we replace $\tilde{\mu}$ with $\mu$. 

The inclusions at $0$ given by $\id_{\boxq}\times \iota_0^{\boxq}$ and $\iota_0^{\boxq}\times \id_{\boxq}$ have image in $\P^1\times \P^1$ that is clearly not disjoint from the center of the blow-up. We explicitly lift them to $B$ by taking the strict transform of their image. Explicitly, for $\iota_0^{\boxq}\times \id_{\boxq}$ (the other case is identical) we have \[(\P^1, \infty) \xrightarrow{\iota_0^{\boxq}\times \id_{\boxq}} (0\times \P^1 \times 0, 0\times \infty\times 0) \hookrightarrow B.\] Then $\partial B \cap  (0\times \infty\times 0) = (0\times \infty\times 0)$, so that the map is admissible. The composition  $(0\times \P^1 \times 0) \hookrightarrow B\xrightarrow{p_3} \P^1$ is the constant morphism to $0\in \P^1$, so that we have identities 
\[\tilde{\mu} \circ (\id_{\boxq}\times \iota_0^{\boxq}) = \tilde{\mu}\circ(\iota_0^{\boxq}\times \id_{\boxq}) = \iota_0^{\boxq} \circ p_{\boxq} \quad \text{ in } \Smlog, \]
that descend to the corresponding identities in $\BSmlog$ once we replace $\tilde{\mu}$ with $\mu$. To summarize, we have proved the following
\begin{prop} The quintuple $(\boxq,  \iota_0^{\boxq },  \iota_1^{\boxq }, p_{\boxq}, \mu)$ makes $\boxq$ into an interval object for the category $\BSmlog$.
    \end{prop}

\subsection{Topologies on $\Smlog$ and $\MSmlog$}\label{subsec:topology-Nis-Zar-modulus-cat} The definitions of this section are adapted from \cite{KSY2} to our setting. Let $\sigma$ be either the Zariski or the Nisnevich topology on $\Sm$.
\begin{df} A morphism $p\colon (U, \partial U)\to (X, \partial X)$ in $\Smlog$ is called a $\sigma$-covering for $(X, \partial X)$ if $p\colon U\to X$ is a $\sigma$-cover of $\Sm$ and   a minimal morphism in $\Smlog$. 
    \end{df}
By Corollary \ref{cor:quarrable_smooth_morphisms}, the pull-back of a $\sigma$-covering along any morphism $f\colon (Y,\partial Y)\to (X, \partial X)$ is still a $\sigma$-covering, so that the above definition gives rise to a Grothendieck topology $t_\sigma$ on $\Smlog$.

The topology $t_\sigma$ is the Grothendieck topology associated with a cd-structure $P_\sigma$. The distinguished squares $P_\sigma$ are defined as follows. Let $(X;\partial X)$ be a pair in $\Smlog$. For $\sigma = {\rm Zar}$,   let $i\colon U\hookrightarrow X$ and $j\colon V\hookrightarrow X$ be two open embeddings (for the Zariski topology on $X$). Then we have the pairs $(U, \partial U)$ and $(V, \partial V)$ in $\Smlog$, where $\partial U = U\cap \partial X$ and $\partial V = V\cap \partial X$ are strict normal crossing divisors on $U$ and $V$ respectively. The distinguished squares $P_{\rm Zar}$ on $\Smlog$ over $(X, \partial X)$ are then given by the pull-back squares
\[\xymatrix{ (U\cap V, \partial X \cap (U\cap V)) \ar[r]\ar[d] & (U, \partial U)\ar[d] \\
(V, \partial V) \ar[r] & (X, \partial X)
}
\]
for $U$ and $V$ running on the set of open subschemes of $X$.

In a similar fashion, an elementary Nisnevich square (i.e. a distinguished square in $P_{\rm Nis}$) in $\Smlog$ is a pull-back square of the form
\[\xymatrix{ (U \times_X Y, \partial Y \cap p^{-1}(U)) \ar[r]\ar[d] & (Y, \partial Y)\ar[d]^p \\
(U, \partial U) \ar[r]^j & (X, \partial X)
}\]
where $j\colon U\hookrightarrow X$ is an open embedding, $p\colon Y\to X$ is an \'etale morphism such that \[(p^{-1}(X\setminus U))_{\rm red} \to (X\setminus U)_{\rm red}\] is an isomorphism, $\partial Y$ is the strict normal crossing divisor $p^{-1}(\partial X)$ on $Y$ and $\partial U = U\cap \partial X$.

The following Proposition is an immediate application of \cite[Lemma 2.5 and Lemma 2.11]{Vcdtop}, using the known results for the Nisnevich and Zariski topology on $\Sm$.
\begin{prop}\label{prop:Nis-Zar-are-regu-comple-cd}The set of elementary Nisnevich (resp. Zariski) squares $P_{\rm Nis}$ (resp. $P_{\rm Zar}$) on the category $\Smlog$ defines a complete and regular cd-structure. In particular, a presheaf of sets $F$ on $\Smlog$ is a sheaf in the Nisnevich (resp. Zariski) topology if and only if for any elementary square $Q$, the square of sets $F(Q)$ is cartesian.
    \end{prop}
We can add the modulus divisor to the picture, obtaining two complete and regular cd-structures on the category of modulus data.
\begin{df}Let $\sigma\in \{{\rm Zar}, {\rm Nis}\}$. A morphism $p\colon U\to M$ of modulus data is called a $\sigma$-cover if 
    \begin{romanlist} 
        \item The underlying morphism of schemes $\ol{p}\colon \ol{U}\to \ol{M}$ is a $\sigma$-cover of $\Sm$,
        \item $p$ is a minimal morphism of modulus data (so that $D_U = p^*(D_M)$ and $\partial U = (p^{-1})(\partial M)$).
        \end{romanlist}
    \end{df}
The class of $\sigma$-covers defines a Grothendieck topology on $\MSmlog$, using Proposition \ref{prop:quarrable_smooth_morphisms} instead of Corollary \ref{cor:quarrable_smooth_morphisms}. The topology $t_\sigma$ on $\MSmlog$ is the Grothendieck topology associated with a complete and regular cd-structure $P_\sigma$. For $\sigma={\rm Nis}$, a distinguished square in $\MSmlog$ is a pull-back square of the form
\begin{equation}\label{eq:Nis-square-MSmlog}\xymatrix{ (\ol{U} \times_{\ol{M}} \ol{Y}, \partial Y \cap p^{-1}({\ol{U}}), D_{\ol{U} \times_{\ol{M}} \ol{Y}}) \ar[r]\ar[d] & ({\ol{Y}}; \partial Y, D_Y)\ar[d]^p\\
(\ol{U}; \partial U, D_U) \ar[r]^j & (\ol{M}; \partial M, D_M)
}
\end{equation}
where $j\colon \ol{U}\hookrightarrow \ol{M}$ is an open embedding, $p\colon \ol{Y}\to \ol{M}$ is an \'etale morphism such that \[(p^{-1}(\ol{M}\setminus \ol{U}))_{\rm red} \to (\ol{M}\setminus \ol{U})_{\rm red}\]  is an isomorphism, and 
\[\partial Y = p^-1(\partial M), \quad D_Y = p^* D_M, \quad \partial U = U\cap \partial M, \quad D_U = D_M \cap U.\]

The cd-structures $P_{\rm Nis}$ and $P_{\rm Zar}$ on $\MSmlog$ are also bounded in the sense of Definition \cite[Definition 2.22]{Vcdtop}. A density structure $D_i(-)$ that works for both cd-structures was introduced by Voevodsky in \cite[2]{Vunst-nis-cdh}. In our context, it takes the following form. Recall that a sequence of points $x_0,\ldots, x_d$ of a topological space $X$ is called \textit{an increasing sequence} of \textit{length} $d$ if $x_i\neq x_{i+1}$ and $x_i\in \ol{\{ x_{i+1}\}}$.
\begin{df} Let $M=(\ol{M}; \partial M, D_M)$ be a modulus datum. Define $D_i(M)$ as the class of morphisms of modulus data $j\colon U\to M$ where $U = (\ol{U};\partial U, D_U)$ is the minimal datum associated to a dense open embedding $j\colon \ol{U}\hookrightarrow \ol{M}$ such that for any $z\in \ol{M}\setminus \ol{U}$ there exists an increasing sequence $z = x_0, x_1, \ldots, x_d$ in $\ol{M}$ of length $d$. 
    \end{df}
\begin{lem}\label{lem:Nis-Zar-are-bounded-cd}The assignment $M\mapsto D_i(M)_{i\geq 0}$ defines a density structure on the category of modulus data, that is compatible with the standard density structure on $\Sm$ defined in \cite[2]{Vunst-nis-cdh}. The cd-structures $P_{\rm Nis}$ and $P_{\rm Zar}$ on $\MSmlog$ are bounded with respect to this density structure.
    \begin{proof} Let $\sigma$ be either the Zariski or the Nisnevich topology. We need to show that the density structure $D_i(M)_{i\geq 0}$ is reducing for the cd-structures $P_{\rm Nis}$ and $P_{\rm Zar}$ on $\MSmlog$. By definition, this property depends only on the small site $M_{\sigma}$ attached to any fixed modulus datum $M$. Let $(\ol{M})_\sigma$ be the usual small $\sigma$-site on the underlying scheme $\ol{M}$. Then the forgetful functor $F\colon \MSmlog\to \Sm$ that sends a modulus datum $N$ to the underlying scheme $\ol{N}$ defines an isomorphism of sites $M_\sigma \xrightarrow{F_M} (\ol{M})_\sigma$ (this was already observed by \cite[Lemma 3.9]{KSY2} in the context of modulus pairs). The claim now follows from  \cite[Proposition 2.10]{Vunst-nis-cdh}.
    \end{proof}
    \end{lem}
    Recall that a point of a site $T$ is a functor $x^*\colon \Shv(T)\to \Set$ which commutes with finite limits and all colimits. A site $T$ has \textit{enough points} if isomorphisms can be tested stalkwise, i.e. if there is a set $x_i^*$ of points such that the induced functor $(x_i^*)\colon \Shv(T)\to \prod_i \Set$ is faithful.
    \begin{ques} Does the site $\MSmlog$ with the Nisnevich or the Zariski topology defined above have enough points?
        \end{ques}
\section{Motivic spaces with modulus}
\subsection{Generalities and first definitions}\label{subsec:Notation-terminologysPsh}
Let $\cS$ be the category of simplicial sets, $\cS = \Delta^{\rm op} \Set$, with simplicial function objects $\cS(-,-)$. We write $\Delta^n$ for the standard $n$-simplex $\hom_\Delta(-, [n])$. If $T$ is a site, we write $\Psh(T)$ for the category of presheaves (of sets) on $T$ and $\sPsh(T) = \Delta^{\rm op}\PSh(T)$ for the category of simplicial objects in $\Psh(T)$. 


For every object $U$ in $T$, we denote by $h_U$ (or simply by $U$ if no confusion arises) the Yoneda functor $h_U(X) = \hom_{T}(X,U)$ considered as discrete simplicial set.


\subsection{Monoidal structures on presheaves categories}\label{subsec:Day-convolution} In this section we present some general material on symmetric monoidal structures for simplicial presheaves. We will then specialize these general results for the construction of a closed symmetric monoidal model structure on the category of motivic spaces with modulus. 


\subsubsection{}Let $\cC$ be a small symmetric monoidal category, with tensor product $\tensor$ and unit $\mathds{1}$. There is a natural extension of the monoidal structure on $\cC$ to a symmetric monoidal structure on $\Psh(\cC)$ via \textit{Day convolution} (see \cite{Day}), that makes the Yoneda functor $h_{\cC}$ strong monoidal. The existence of the monoidal structure follows formally from the general theory of left Kan extensions. 
Explicitly, given that any presheaf is colimit of representable presheaves, write $F= \colim_{X \downarrow F}h_X$ and $G = \colim_{Y \downarrow G}h_Y$. Then their convolution $F\Daytensor G$ is the colimit $\colim_{X,Y}h_{X\tensor Y}$. It follows immediately from the definition that the Yoneda functor is strong symmetric monoidal.  It is also clear that the unit for the convolution product is the representable presheaf $h_{\mathds{1}}$. It is also formal to see that the monoidal structure on presheaves given by Day convolution is closed, i.e. there exists an internal hom $[-,-]^{\rm Day}$ that is right adjoint to $\Daytensor$. This is characterized by
\[[F, G]^{\rm Day}(X) = \hom_{\Psh(\cC)}(h_X\Daytensor F, G), \quad  \text{ for all } X\in \cC.\]

Unless required for clarity, we will drop the superscript and write simply $\tensor$ for the tensor product of presheaves.
Recall (see e.g. \cite{ImKelly}) that a monoidal category $(\cD, \star, \mathds{1}_\cD)$ is called monoidally co-complete if $\cD$ is co-complete and all the endofunctors $X\star (-)$, $(-)\star Y$ for $X, Y$ in $\cD$ are co-continuous. By \cite[Proposition 4.1]{ImKelly}, the category $\Psh(\cC)$ on a small symmetric monoidal category is monoidally co-complete. This construction is universal in the following sense.
\begin{prop}[Theorem 5.1 \cite{ImKelly}]\label{prop:ImKelly} Let $\cD$ be a monoidally co-complete category. Then the functor $[\Psh(\cC), \cD]_\tensor \to [\cC, \cD]_\tensor$ between the categories of strong monoidal functors from $\Psh(\cC)$ to $\cD$ and the category of strong monoidal functors from $\cC$ to $\cD$ induced by the Yoneda functor is an equivalence.
    \end{prop}
\begin{rmk}Here's a situation where the previous Proposition turns out to be useful. Let $u\colon \cC \to \cD$ be a strong symmetric monoidal functor. Then $u$ gives rise to a string of adjoint functors between the categories of presheaves   
\[(u_!, u^*, u_*),\quad \Psh(\cD)\xrightarrow{u^*} \Psh(\cC)\]
where each functor is the left adjoint to the the following one. The left adjoint $u_!$ to the restriction $u^*$ is defined via left Kan extension, so that Proposition \ref{prop:ImKelly} implies that it is strong monoidal. 
\end{rmk}
\subsubsection{}There is a natural way of extending Day convolution from the category of presheaves on $\cC$ to the category of simplicial presheaves, so that the sequence of embeddings 
\[\cC\hookrightarrow \Psh(\cC)\hookrightarrow \sPsh(\cC)\]
is a sequence of strong monoidal functors (recall here that we identify presheaves of sets with discrete simplicial presheaves, i.e. simplicial presheaves of simplicial dimension zero). Given two simplicial presheaves $F,G$ we define $F\tensor G$ by 
\begin{equation}\label{eq:extension-Day-psh-to-spsh}(F\tensor G)_n = F_n\tensor G_n, \quad n\geq 0\end{equation}
and this gives $\sPsh(\cC)$ the structure of a closed symmetric monoidal category.  We keep writing $[-,-]$ for the internal hom for Day convolution.
\begin{rmk}\label{rmk:enrichment-simpl-preshaves-Day-and-cartesian-prod}The category of simplicial presheaves on $\cC$ is tensored over $\cS$ in the following way. The product $F\times K$ of a simplicial presheaf $F$ with a simplicial set  $K$ is defined on sections by $(F\times K)(U) = F(U)\times K$. Alternatively, we can simply think to $X\times K$ as the product of $X$ with the constant simplicial presheaf $K$.
    
The functor $\cS\to \sPsh(\cC)$ given by $K\mapsto \mathds{1}\Daytensor K$ is easily seen to be endowed with the structure of strong symmetric monoidal functor.
\end{rmk}
\begin{rmk}\label{rmk:pointed-Day}There is a pointed variant of Day convolution. Let $\sPsh(\cC)_\bullet$ be the category of pointed simplicial presheaves, i.e. the category of presheaves of pointed simplicial sets on $\cC$, and let \[(-)_+\colon \sPsh(\cC)\to \sPsh(\cC)_\bullet\] be the canonical ``add base point'' functor (left adjoint to the forgetful functor). By mimicking the definition of smash product $\wedge$ of pointed simplicial preshaves starting from the cartesian product, we can define a symmetric monoidal structure $\Daytensor_\bullet$ on  $\sPsh(\cC)_\bullet$ (see Section \ref{sec:pointed-spaces} for details). This is the unique symmetric monoidal structure on  $\sPsh(\cC)_\bullet$ that has $\mathds{1}_+$ as unit and that makes $(-)_+$ strong monoidal. 
\end{rmk}
We will come back later on the behaviour of Day convolution with respect to different model structures on simplicial presheaves.
\subsection{Motivic spaces with modulus and interval objects}\label{sec:Mot-spaces-and-intervals}We begin with the following definition.  
\begin{df}\label{def:motivic-spaces}Let $\MSmlog$ be the category of modulus data over $k$. A \textit{motivic space with modulus} is a contravariant functor $\kX\colon \MSmlog\to\cS$ i.e. a simplicial presheaf on $\MSmlog$. We let $\MM$ denote the category of motivic spaces with modulus.
    \end{df}
Since $\MM$ is a category of simplicial presheaves on a small category, it is a locally finitely presentable bicomplete $\cS$-category, with simplicial function complexes defined as above. In particular, finite limits commute with filtered colimits. The following fact is standard. 

\begin{lem}Every motivic space with modulus is filtered colimit of finite colimits of spaces of the form $h_M\times \Delta^n$, for $M\in \MSmlog$ a modulus datum and $n\geq 0$.
    \end{lem}
Apart from the category of motivic spaces $\MM$, there are two other categories of simplicial presheaves that will play an important r\^ole in what follows.
\begin{df}The category of \textit{motivic spaces with compactifications}, $\Mlog$ is the category of simplicial presheaves on $\Smlog$, i.e. $\Mlog = \sPsh(\Smlog)$. The category of \textit{birational motivic spaces with compactifications}, $\BMlog$ is the category of simplicial presheaves on the localized category $\BSmlog$, i.e. $\BMlog = \sPsh(\BSmlog)$. Finally, we let $\cM_k$ denote the category of \textit{motivic spaces} over $k$ in the sense of Morel-Voevodsky, i.e. $\cM_k = \sPsh(\Sm)$.
\end{df}
The categories $\MM$, $\Mlog$ and $\BMlog$ are closed symmetric monoidal categories, where we consider on $\MM$ Day convolution induced by the monoidal structure \ref{subsec:monoidal-structure-MSmlog} on $\MSmlog$ (and the usual Cartesian product on the other categories).
\subsubsection{}Recall from the discussion in Section \ref{subsec:inv-rat-maps} (with the notations of \eqref{eq:loc-functorSmlog}) that there is a canonical faithful functor \[v\colon \Smlog\to \BSmlog\] and from \ref{subsubsec:vfunctor} that there is a fully faithful embedding \[u\colon \Smlog\hookrightarrow \MSmlog.\] They are both strict monoidal functors.

These functors extend to the presheaves categories, giving a plethora of adjunctions
\[(u_!, u^*, u_*),\quad u^*\colon \MM\leftrightarrows \Mlog\colon u_* \]
\[(v_!, v^*, v_*),\quad v^*\colon \BMlog \leftrightarrows\Mlog\colon v_*.\]
Note that from  general principle the restriction functors $u^*$ and $v^*$ preserve  limits and colimits, and the functors $u_!$ and $v_!$  preserve colimits  and are strong monoidal by Proposition \ref{prop:ImKelly}. 
\begin{df}\label{def:interval-I}We denote by $\cI$ the object of $\MM$ given by 
    \[\cI= u_! v^* (\boxq) = u_! v^* ( h_{(\P^1,\infty)}).\]
    \end{df}
We will use the interval structure of $\boxq = \boxq^1$ on $\BSmlog$ to show that $\cI$ is an interval object in the symmetric monoidal category $\MM$. We start from the following simple observation.
\begin{lem}The representable simplicial presheaf $\boxq$ is an interval object in $\BMlog$ and a weak interval object in $\MM$.
    \begin{proof}The maps $(\iota^\boxq_0, \iota^\boxq_1, p_\boxq)$ define maps in $\BMlog$. Since the Yoneda embedding preserves (small) limits, we have $h_\boxq \times h_\boxq = h_{\boxq\times \boxq} = h_{\boxq\tensor \boxq}$, so that the multiplication $\mu$ also extends, with the required compatibilities, to $\BMlog$. The statement for $\MM$ is also clear.
        \end{proof} 
    \end{lem}

\subsubsection{}\label{subsubsec:ComparisonMapIntervals}Let $(X, \partial X)$ be an object of $\Smlog$. From the adjunction $(v_!, v^*)$ we get a natural map $\eta_X\colon h_{(X, \partial X)}\to v^*(v_! h_{(X, \partial X)}) = v^*( h_{(X,\partial X)})$ (since $v_!$ commutes with Yoneda). Evaluated on an object $(Y, \partial Y)$ of $\Smlog$, the map $\eta_X$ corresponds to the inclusion
\[ \hom_{\Smlog}((Y, \partial Y), (X, \partial X))\hookrightarrow \hom_{\BSmlog}((Y, \partial Y), (X, \partial X)).\]
We will still denote by $\eta_X$ the morphism in $\MM$ given by $u_!(\eta_X)$. This is the map of motivic spaces 
\[ h_{(X; \partial X, \emptyset)} = u_!(h_{(X, \partial X)})\xrightarrow{u_!\eta_X} u_!(v^*(h_{(X, \partial X)})).\]
For $X = \boxq$, the above construction gives a canonical \textit{comparison} morphism of motivic spaces $\eta\colon \boxq \to \cI$.
\begin{prop}\label{prop:I-is-interval}The motivic space $\cI$ is an interval object in $\MM$ for the Day convolution product.
    \begin{proof} We start by proving that the simplicial presheaf $v^*(\boxq)$ in $\Mlog$ is an interval object for the usual product of presheaves. Since the terminal object of $\Smlog$ is $(\Spec{k}, \emptyset)$ and since Yoneda preserves small limits, the simplicial presheaf represented by $(\Spec{k}, \emptyset)$ is just the constant simplicial set having one element in simplicial degree $0$ (the ``point'', denoted $\pt$). Since $v^*$ preserves finite limits, we have that $v^*(\pt) = \pt$, and therefore we automatically obtain maps
        \[\iota_\epsilon^{v^*(\boxq)}\colon \pt = v^*(\pt)\rightrightarrows v^*(\boxq), \quad \epsilon\in \{0,1\}, \quad p_{v^*(\boxq)}\colon v^*(\boxq)\to \pt\]
        satisfying the identities $p_{v^*(\boxq)} \circ \iota_\epsilon^{v^*(\boxq)} = \id_{\pt}$. Let now $\mu\colon \boxq\times \boxq\to \boxq$ be the multiplication map in $\BMlog$. Since again $v^*$ preserves finite limits, we have $v^*(\boxq\times \boxq) = v^*(\boxq)\times v^*(\boxq)$ so that we get a map
        \[v^*\mu \colon v^*(\boxq)\times v^*(\boxq)\to v^*(\boxq)\]
        and we have
        \begin{align*} v^*(\mu)\circ (\id_{v^*\boxq} \times v^*(\iota_0^{\boxq})) =& v^*(\mu)\circ v^*(\id_{\boxq} \times (\iota_0^{\boxq})) = v^*(\mu\circ (\id_{\boxq} \times (\iota_0^{\boxq})) \\=& v^*(\iota_0^\boxq\circ p_\boxq) = \iota_0^{v^*\boxq} \circ p_{v^*\boxq},
        \end{align*}
        \begin{align*} v^*(\mu)\circ (\id_{v^*\boxq} \times v^*(\iota_1^{\boxq})) =& v^*(\mu)\circ v^*(\id_{\boxq} \times (\iota_1^{\boxq})) = v^*(\mu\circ (\id_{\boxq} \times (\iota_1^{\boxq})) =\id_{v^*{\boxq}},
        \end{align*}
    completing the proof that $v^*(\boxq)$ is an interval in $\Mlog$. As for $\cI$, we first notice that \[u_!(h_{(\Spec{k}, \emptyset)}) = \mathds{1}\] is the unit for Day convolution on $\MM$. Applying the functor $u_!$ we then obtain the morphisms $\iota_0^{\cI}, \iota_1^{\cI}, p_{\cI}$ that make $\cI$ into a weak interval on $\MM$. As for the multiplication, it's enough to show that $u_!(v^*(\boxq) \times v^*(\boxq)) = u_!(v^*(\boxq\times \boxq)) \cong \cI\tensor \cI$, since the identities involving $u_!(\mu)$ will be then automatically satisfied by functoriality (or using the same chain of equalities as above). Slightly more generally, let $F\in \Psh(\Smlog)$ be a presheaf of sets and consider it as simplicial presheaf of simplicial dimension zero. Write $F = \colim_{U\downarrow F}h_U$, for $U\in \Smlog$. Then
    \begin{align*} u_!(F\times F)  &= u_!(\varvarcolim_{U\downarrow F} h_U\times\varvarcolim_{U'\downarrow F} h_{U'}) =^{\dagger} u_!(\varvarcolim_{U,U'} h_{U\times U'}) \\ {}^{\dagger \dagger}&= \varvarcolim_{U, U'} u_!(h_{U\times U'}) = \varvarcolim_{U, U'} (h_{u_!(U\times U')}) =  \varvarcolim_{U, U'} (h_{u(U)\tensor u(U')})  \\ 
    & = \varvarcolim_{U,U'} (h_{u(U)}\tensor h_{u(U')}) =^{\ddagger} \varvarcolim_{U\downarrow F} h_{u(U)}\tensor \varvarcolim_{U'\downarrow F} h_{u(U')}\\ &= u_!(F)\tensor u_!(F).
  \end{align*}
 The equality $\dagger$ follows from the fact that colimits commute with finite fiber products in a category of presheaves (Giraud's axiom), while $\dagger\dagger$ follows from the fact that $u_!$ commutes with colimits. For the equality $\ddagger$ we have used the fact that $\Psh(\MSmlog)$ is monoidally co-complete for Day convolution product. The other equalities are trivial, using the fact that $u_!$ commutes with Yoneda and that $u(U)\tensor u(U') = u(U)\times u(U')$ in $\MSmlog$ for every $U, U'$ in $\Smlog$ by \ref{subsubsec:tensor-cart-empty-modulus}. Specializing these equalities to the case $F= v^*(\boxq)$ gives the required statement.
 \end{proof}
    \end{prop}
    \begin{rmk} Our method of transporting the interval structure from $\BSmlog$ to $\MM$ looks quite general. It seems plausible that one can repeat a similar argument by replacing the category of simplicial presheaves on $\BSmlog$ with the category of extended co-cubical objects in $\cS$, i.e. the category of strong monoidal functors $[\mathbf{ECube}, \cS]_\tensor$, or even with the category of extended cubical object in a category of presheaves with values in monoidal model category $\cM$. 
\end{rmk}
\begin{rmk}In the references \cite{MV}, \cite{VSelecta} and \cite{Blander} (for $\A^1$-theory) and \cite{KSY2} (for the modulus-theory), the interval objects considered are always representable (either by $\A^1$ or by the modulus pair $(\P^1, 1)$). Here, we are pushing the ideas of \cite{KSY2} further to get a theory that works more generally for interval objects in categories of presheaves that are not necessarily representable.
    \end{rmk}
\subsubsection{}\label{subsubsec:canonical-maps-btw-intervals}Let $\A^1$ denote the representable simplicial presheaf $h_{(\A^1;\emptyset, \emptyset)} = h_{u(\A^1, \emptyset)}$. It is clearly an interval object in $\MM$ for the cartesian product as well as for Day convolution product, since by \ref{subsubsec:tensor-cart-empty-modulus} we have $\A^1\tensor \A^1 = \A^1\times \A^1$.  
From the admissible morphism $j\colon (\A^1, \emptyset)\hookrightarrow (\P^1, \infty)$ in $\Smlog$ we get maps of motivic spaces \[\A^1\to \boxq\to \cI \text{ in }\MM,\quad \A^1\to \boxq\to v^*(\boxq) \text{ in }\Mlog \]
that we will use to compare the different interval structures on $\MM$, on $\Mlog$ and on $\cM(k)$.

\subsubsection{}We start by comparing $\A^1$ and $v^*(\boxq)$ with the standard interval $\A^1$ in the category of motivic spaces $\cM(k)$. Recall that the adjoint pair $\lambda\colon \Sm\leftrightarrows \Smlog\colon \omega$ gives rise to a string of four adjoint functors
\[(\lambda_!, \lambda^* = \omega_!,\lambda_* = \omega^*, \omega_*), \quad \omega_!=\lambda^*\colon \Mlog\leftrightarrows \cM(k)\colon \lambda_*\]
(i.e. the functor $\omega_!$ has in turn a left adjoint). The functors $\lambda^*$ and $\omega^*$ clearly commute with products. Since $\lambda$ commutes with products and $\lambda_!$ commutes with Yoneda, the same argument used in the proof of Proposition \ref{prop:I-is-interval} (or even Proposition \ref{prop:ImKelly} directly) shows that $\lambda_!$ is monoidal with respect to the cartesian product. 
\begin{lem}\label{lem:comparison-A1-box-M}There are canonical isomorphisms $\lambda^*(\A^1) \cong \lambda^* v^*(\boxq) \cong \A^1$ as interval objects of $\cM(k)$.
    \begin{proof}Since $\lambda^* = \omega_!$, it's clear that $\lambda^*(\A^1) = \omega_!(h_{(\A^1, \emptyset)}) = h_{\omega(\A^1, \emptyset)} = h_{\A^1} = \A^1$. For the second statement, it's enough to check that $\lambda^*v^*(\boxq) \cong\A^1$ as presheaves on $\Sm$, and we just have to play with adjunctions. For $X\in \Sm$, we have
        \begin{align*}\hom_{\Psh(\Sm)}( h_X, \lambda^*v^*(\boxq)) &= \hom_{\Psh(\Smlog)}(\lambda_!(h_X), v^*(\boxq)) \\ &= \hom_{\Psh(\BSmlog)}(v_!(h_{(X, \emptyset)}), \boxq) \\ &=\hom_{{\BSmlog}}((X, \emptyset), (\P^1, \infty)) \\ &= \hom_{\Sm}(X,  \A^1).
        \end{align*}
The fact that the isomorphisms are compatible with the interval structure is a tautology from the definitions.
        \end{proof}
    \end{lem}
\begin{rmk}Note that the interval $\cI$ of $\MM$ is obtained from $v^*(\boxq)$ by applying the left adjoint $u_!$ of the restriction functor $u^*$. The exact same construction, using $\omega_!=\lambda^*$ instead of $u_!$ produces, in view of Lemma \ref{lem:comparison-A1-box-M}, nothing but the usual interval $\A^1$ in the Morel-Voevodsky category of motivic spaces.
    \end{rmk}
\subsubsection{}Before moving forward to give more refined comparisons, we ask the following question. Let $F\colon \Smlog\to \Sm$ be the forgetful functor $(X,\partial X)\to X$. Then $F$ defines the usual set of adjoint functors between the categories of spaces $\Mlog$ and $\cM(k)$, namely $(F_!, F^*, F_*)$. The general principle illustrated above allows us to construct the object $J\in \cM$ as $J=F_!(v^*(\boxq))$. Since $F$ is strict monoidal for the cartesian product, the argument of Proposition \ref{prop:I-is-interval} goes through to show that $J$ is in fact an interval object in $\cM(k)$ for the usual product of simplicial presheaves.

The canonical adjunction map $\eta\colon \boxq\to v^*(\boxq)$ gives then a map $\P^1\to J$ in $\cM$. By construction, this map cannot be an isomorphism (since $\P^1$ does not have any interval structure on $\cM(k)$), and is also easy to see that $J$ is not isomorphic to $\A^1$. One could therefore try to develop a machinery for the localization of $\cM(k)$ to $J$, in the same spirit of what we will do for $\MM$ with $\cI$. We don't know, at the moment, what would be the outcome of such a construction. 

\section{Motivic homotopy categories with modulus}
\subsection{Model structures on simplicial presheaves}
The category of simplicial sets $\cS$ carries a well-known cofibrantly generated model structure (the Quillen model structure). 
The category of simplicial presheaves $\sPsh(T)$ on a small Grothendieck site $T$ carries two natural model structures, the \textit{injective} and the \textit{projective} model structure presenting the same homotopy category, i.e. having the same class of weak equivalences:
\begin{df}A map $f\colon A\to B$ in $\sPsh(T)$ in $\sPsh(T)$ is called an \textit{objectwise} (or \textit{levelwise} or \textit{sectionwise}) \textit{simplicial weak equivalence} if $f(X)\colon A(X)\to B(X)$ is a weak equivalence of simplicial sets for each $X\in T$. A map $f\colon A\to B$ is called an \textit{objectwise Kan fibration} if $f(X)\colon A(X)\to B(X)$ is a Kan fibration of simplicial sets for each $X\in T$.
    \end{df}
 The injective model structure on $\sPsh(T)$ was introduced by Heller in \cite{Heller}, in which the cofibrations are the monomorphisms, and the fibrations are the maps having the right lifting property with respect to the trivial cofibrations. It is a proper, simplicial cellular model structure.
We will refer to the fibrations for the injective model structure as the \textit{(simplicial) injective fibrations}. The category of simplicial presheaves with the injective model stucture will be denoted, following the usual convention, $\sPsh(T)_{\rm inj}$. Note that every object is cofibrant for the injective structure.
The second model structure on $\sPsh(T)$, the projective one, goes back to Quillen (see \cite{Hirsch}, Theorem 11.6.1). In this model structure, the fibrations are the objectwise Kan fibrations, and the cofibrations are the maps having the left lifting property with respect to the trivial fibrations. The projective model structure is also proper, simplicial and cellular. 
We denote by $\sPsh(T)_{\rm proj}$ the category of simplicial presheaves with the projective model structure.  We refer the reader to \cite[12]{Hirsch} for the definition of a cellular model category. Both the injective and the projective model structure on simplicial presheaves are cellular (see \cite{JardineLocal}, around 7.19 for a comment on the cellularity of the injective model structure). By definition,  a cellular model category is cofibrantly generated. 
We choose a functorial cofibrant replacement $(-)^{\rm c}\to \id_{\sPsh{(T)}}$ for the projective model structure, so that for every  object $\kX \in \sPsh(T)$, there is an objectwise trivial fibration $\kXc\to \kX$ with $\kXc$ cofibrant. 
     \subsubsection{}Recall (see \cite[Definition 4.2.6]{Hovey}) that a model category $\cM$ that is also a monoidal category with product $\tensor$ and unit $\mathds{1}$ is called a \textit{monoidal model category} if the following conditions are satisfied
     \begin{enumerate}
         \item Let $Q\mathds{1}\xrightarrow{q}\mathds{1}$ be the cofibrant replacement for the unit $\mathds{1}$. Then the natural map $Q\mathds{1}\tensor X\to \mathds{1}\tensor X\cong X$ is a weak equivalence for all cofibrant $X$.
         \item Given two cofibrations $f\colon U\to V$ and $g\colon W\to X$ in $\cM$, their push-out product
         \[f\square g\colon (V\tensor W) \amalg_{U\tensor W} (U\tensor X)\to V\tensor X\]
         is a cofibration, which is trivial if either $f$ or $g$ is.
         \end{enumerate}
         We refer to the first condition as the \textit{unit axiom} and to the second condition as the \textit{pushout product axiom}.
                 
Suppose that a small site $T$ carries a symmetric monoidal structure $\tensor$, that extends via  Day convolution to $\sPsh(T)$ (see Section \ref{subsec:Day-convolution}). The proof of the following proposition is easy, using \cite[Corollary 4.2.5]{Hovey}.
\begin{prop}\label{prop:Day-is-monoidal} The projective model structure on simplicial presheaves is a (symmetric) monoidal model category with respect to Day convolution.
\end{prop}
\begin{rmk}The pointed version of Day convolution presented in \ref{rmk:pointed-Day} gives rise to a symmetric monoidal model structure on the category of pointed simplicial presheaves $\sPsh(T)_\bullet$ with the projective model structure. This follows from Proposition \ref{prop:Day-is-monoidal} together with \cite[Proposition 4.2.9]{Hovey}. Alternatively, one can repeat the proof of Proposition \ref{prop:Day-is-monoidal} replacing the product of simplicial presheaves with the smash product of the pointed counterparts.
    \end{rmk}
    
\subsection{Local model structures} Following Jardine (see \cite{JardineLocal} or \cite{JardineJPAA}), we can put the topology in the picture as follows. Let $T$ be again a small Grothendieck site and let $\sPsh(T)$ be the category of simplicial presheaves on $T$. 
The reader can consult \cite{JardineLocal}, Section 4.1 and Section 5.1 for the notion of \textit{local weak equivalence} and for the notion of \textit{injective} or \textit{global} fibration associated to the topology $\sigma$ on $T$ for a morphism $f\colon \kX\to \kY$ in $\sPsh(T)$.  Recall finally that an \textit{injective cofibration} is simply a monomorphism of simplicial presheaves. If necessary, we will distinguish injective fibrations for the injective-local model structure from fibrations for the injective model structure on simplicial presheaves by adding the word \textit{simplicial} to the latter class. By \cite{JardineLocal}, Theorem 5.8, the local injective model structure $\sPsh(T)_{\rm inj}^{loc}$  (having local weak equivalences as weak equivalences and injective fibrations as fibrations) is a proper simplicial cellular closed model category. 
\subsubsection{Localizations}Given any model category $\cM$ and a set of morphisms $S$, we say that an object $X$ of $\cM$ is $S$-local (see \cite[3.1.4]{Hirsch}) if it is fibrant and for every map $f\colon A\to B$ the induced map $f^*$ between the homotopy function complexes is a weak equivalence. A map $g\colon X\to Y$ is called an $S-$local equivalence if for every $S$-local object $Z$ the induced map $g^*$ between the homotopy function complexes is a weak equivalence.  We write $L_{S}\cM$  for the model structure on $\cM$ in which the weak equivalences are the $S$-local equivalences of $\cM$, the cofibrations are the same cofibrations of $\cM$ and the class of fibrations is the class of maps having the right lifting property with respect to those maps that are both cofibrations and $S$-local equivalences. This is called the \textit{left Bousfield localization of $\cM$  with respect to $S$}. By a general result of Hirschhorn, see \cite[Theorem 4.1.1]{Hirsch}, if $\cM$ is left proper and cellular, $L_S{\cM}$ exists and is a left proper cellular model category, that has a natural  structure of simplicial model category if $\cM$ has one. 
        
\subsubsection{}Let $\cM$ be a (symmetric) monoidal model category. Write $\tensor$ for its monoidal product and take a set of morphisms $S$. In general, the Bousfield localization $L_{S}\cM$ (whenever exists) will not inherit the structure of monoidal model category. The following Proposition gives a convenient criterion for checking if the localization to $S$ behaves well with respect to the monoidal structure. The proof is standard, and we refer to \cite{CT}.
\begin{prop}[\cite{CT}, Proposition 5.6]\label{prop:loc-monoidal-structure-CT}Let $\cM$ be  a left proper cellular symmetric monoidal model category. Let $S$ be a set of morphisms in $\cM$. Assume that the following conditions hold:
    \begin{romanlist}
        \item $\cM$ admits generating sets of (trivial) cofibrations  consisting of maps between cofibrant objects;
        \item For every cofibrant object $X$, the functor $X\tensor(-)$ sends the elements of $S$ to $S$-local weak equivalences;
        \item The unit object for the monoidal structure is cofibrant.
        \end{romanlist}
        Then the left Bousfield localization $L_S\cM$ with respect to $S$ (that exists by  \cite[Theorem 4.1.1]{Hirsch}) is a symmetric monoidal model category.
    \end{prop}
\subsubsection{}\label{sec:def-SigmaP-def-proj-is-bousfield}When the topology on $T$ is defined by a regular, complete and bounded cd-structure $P$, a result of Voevodsky (\cite{Vcdtop}, Proposition 3.8) presents the local injective structure of Jardine as left Bousfield localization of the injective structure $\sPsh(T)_{\rm inj}$ to the class of maps given by distinguished squares. More precisely,  let $T$ be a site as above and let $X\in T$. For $Q$ a distinguished square of the form 
        \begin{equation}\label{eq:squarecdtop}\xymatrix{ B \ar[d] \ar[r]^{e_B} & Y\ar[d]^p\\
         A\ar[r]^e & X
         }
        \end{equation}
 for the cd-structure $P$ on $T$, write $P(Q)$ for the simplicial homotopy push-out of the diagram $(A\leftarrow B\rightarrow Y)$. There is a natural map $P(Q)\to X$. 
Write $\Sigma_P$ for the class of maps
\[\Sigma_P = \{P(Q)\to X\}_{Q}\cup \{\emptyset \to h_\emptyset\}\]
where $Q$ runs on the set of distinguished squares for the cd-structure $P$ on $T$, and $\emptyset$  is the initial object of $\Psh(T)$ and $h_\emptyset$ is presheaf represented by the initial object of $T$.
\begin{thm}[\cite{Vcdtop}]\label{thm:inj-loc-structure-is-bousfield}Let $T$ be a small site whose topology is defined by a complete bounded and regular cd-structure. Then the local injective model structure on $\sPsh(T)$ is the left Bousfield localization of the (global) injective model structure $\sPsh(T)_{\rm inj}$ to the class $\Sigma_P$.
    \end{thm}
\subsubsection{}Together with the injective local model structure, there is a \textit{projective local model structure} on simplicial presheaves due to Blander \cite{Blander}. We recall here the following useful facts about it.
\begin{prop}[\cite{Blander}, Lemma 4.1]\label{prop:BlanderProp} Let $T$ be a site with an initial object $\emptyset$ whose topology is defined by a complete bounded regular cd-structure $P$. Then a simplicial presheaf $F$ on $T$ is local projective fibrant if and only if $F(U)$ is a Kan simplicial set for all $U$ in $T$ and if for every distinguished  square $Q$ of $P$, the square $F(Q)$ is a homotopy pull-back square of simplicial sets. Such presheaves are called \textit{flasque}.
    \end{prop}
    \begin{thm}[Blander \cite{Blander}, Lemma 4.3, Voevodsky \cite{Vcdtop}]\label{thm:Blander-proj-is-Bousfield} The local projective model structure on the category of simplicial presheaves on a site $T$ equipped with a complete regular bounded cd-structure $P$ is the left Bousfield localization of the projective model structure  $\sPsh(T)_{\rm proj}$ to the class $\Sigma_P$.
        \end{thm}
        We write $\sPsh(T)_{\rm proj}^{loc}$ to denote the local projective model structure on $\sPsh(T)$. By the general theory of Bousfield localization, the local projective and the local injective model structures on $\sPsh(T)$ are both cellular, proper and simplicial. 
    
    The weak equivalences in both model structure agree and are precisely the local weak equivalences. Note that left proper is automatic in both cases by Bousfield localization, while right properness follows from \cite[Lemma 3.4]{Blander} (for the projective case) and \cite[Lemma 4.37]{JardineLocal} (for the injective case). 

\begin{prop}\label{prop:Day-makes-proj-loc-monoidal}Let $T$ be a small site whose topology is defined by a complete bounded and regular cd-structure. Assume that for every distinguished square $Q$ in $T$ and every object $Z$ in $T$, the product $Q\tensor Z$ is still a distinguished square in $T$. Then Day convolution makes $\sPsh(T)_{\rm proj}^{loc}$ a monoidal model category.
    \begin{proof}
      Thanks to Blander's Theorem \ref{thm:Blander-proj-is-Bousfield}, the local projective model structure is a (left) Bousfield localization of the projective structure, and we have an explicit description of the set of maps that we are inverting. In the notations of \ref{sec:def-SigmaP-def-proj-is-bousfield}, write $\varphi_X\colon P(Q)\to X$ for the natural map from the simplicial homotopy push-out of a distinguished square having $X\in T$ as bottom right corner. According to \cite[Theorem 5.7]{CT}, we have to check that for every representable presheaf $h_Z$, the induced map $\varphi_X\tensor \id_Z\colon P(Q)\tensor h_Z\to h_X\tensor h_Z = h_{X\tensor Z}$ is an $S$-local equivalence for $S=\Sigma_P$. Write $Q_Z$ for the tensor product square $Q\tensor Z$ in $T$. By assumption, $Q_Z$ is a distinguished square, and the bottom right corner of it is the product $X\tensor Z$. The natural map $P(Q_Z)\to h_{X\tensor Z}$ is an element of $\Sigma_P$ by assumption, and factors through the product $P(Q)\tensor h_Z$. Since the projective model structure is monoidal for Day convolution, the tensor product $(-)\tensor h_Z$ commutes with homotopy colimits. In particular, the map $P(Q_Z)\to P(Q)\tensor h_Z$ is a weak equivalence. By the 2 out of 3 property of $S$-local equivalences, we conclude that the required map $\varphi_Z\tensor \id_Z$ is an $S$-local equivalence, completing the proof.
    \end{proof}
    \end{prop}
\subsection{Interval-local objects and $I$-homotopies} We start from the following general definition.  Let $T$ be a (small) site and let $\sPsh(T)$ be the category of simplicial presheaves on $T$. Suppose that $T$ carries a symmetric monoidal structure $\tensor$, that extends via Day convolution to $\sPsh(T)$. Finally, suppose that there exists a weak interval object $I$ for the $\tensor$-structure on $\sPsh(T)$. We consider on $\sPsh(T)$ both the injective and the projective local model structures.
\begin{df}\label{def:I-local-proj-inj}A simplicial presheaf $\kX$ is called \textit{projective $I$-local} (resp. \textit{injective $I$-local}) if:
    \begin{romanlist} \item The presheaf $\kX$ is fibrant for the projective local model structure on $\sPsh(T)$ (resp. $\kX$ is fibrant for the injective local model structure on $\sPsh(T)$), and
        \item For every $\kY$ in $\sPsh{(T)}$, the map between the homotopy function complexes
        \[(\id\tensor \iota_0^I)^*\colon \Map(\kY\tensor I, \kX)\to \Map(\kY, \kX)\]
induced by $\id\tensor \iota_0^I$ is a weak equivalence.
    \end{romanlist}
    \end{df}
Note that, since $I$ is a weak interval, we have the identities $p_I\circ \iota_0^I=p_I\circ \iota_1^I = \id$ in $T$. In particular, for every $\kY\in \sPsh{(T)}$ the composition
\[\Map(\kY, \kX) \xrightarrow{(\id \tensor p_I)^*} \Map(\kY\tensor I, \kX)\xrightarrow{(\id\tensor \iota_0^I)^*} \Map(\kY, \kX) \]
 is the identity. If $\kX$ is projective $I$-local, then the second map is also a weak equivalence by condition ii) above, and therefore~---~by the 2 out of 3 property of the class of simplicial weak equivalences~---~so is the map $(\id \tensor p_I)^*$. Thus one can replace condition ii) with the following
 \begin{itemize}
     \item[ii')] For every $\kY$ in $\sPsh{(T)}$, the map between the homotopy function complexes
     \[(\id\tensor p_I)^*\colon \Map(\kY, \kX) \to \Map(\kY\tensor I, \kX)\]
     induced by $\id\tensor p_I$ is a weak equivalence.
     \end{itemize}
In a similar way, one can replace condition ii) with the analogue condition stated using the map $\id\tensor \iota_1^I$ instead. All three notions are equivalent. 
\begin{rmk}If $\cM$ is a simplicial model category, a homotopy function complex between a cofibrant object $X$ and a fibrant object $Y$ is weakly equivalent to the simplicial  mapping space (or simplicial function complex, in the terminology of \ref{subsec:Notation-terminologysPsh}) $\cS(X,Y)$. Since every object is cofibrant in $\sPsh(T)_{\rm inj}$, an injective fibrant simplicial presheaf $\kX$ is injective $I$-local if for every $\kY$, the natural map between the simplicial function complexes
    \[\cS(\kY\tensor I, \kX)\to \cS(\kY, \kX)\] 
  is a weak equivalence.
  
  On the other hand, not every object is cofibrant for the projective structure. We can then reformulate condition ii) of Definition \ref{def:I-local-proj-inj} using the simplicial function complex and the functorial cofibrant replacement $(-)^{\rm c}$ as follows: a projective fibrant simplicial presheaf $\kX$ is projective $I$-local if for every object $\kY$, the natural map between the simplicial function complexes
    \[\cS((\kY)^{\rm c}\tensor I^{\rm c}, \kX)\to \cS((\kY)^{\rm c}, \kX)\]
    is a weak equivalence. Note that we are using here the fact that $\tensor$ is left Quillen bi-functor, so that it preserves cofibrant objects. 
    \end{rmk}
\begin{rmk} Note that the cofibrant replacement $(I)^{\rm c}\to I$ can be chosen to be a monoidal functor. This gives automatically the object $(I)^{\rm c}$ the structure of interval object in $\sPsh(T)_{\rm proj}$. This applies, in particular, in the case $T= \MSmlog$ and $\sPsh(T) = \MM$.
    \end{rmk}
\begin{df}\label{def:I-we} A morphism $g\colon \kX\to \kY$ is called a \textit{projective $I-$weak equivalence} (resp.~\textit{an injective $I$-weak equivalence}) if for any projective (resp. injective) $I$-local object $\kZ$ the induced map between the homotopy function complexes
    \[g^*\colon \Map(\kY, \kZ) \to \Map(\kX, \kZ)\]
    is a weak equivalence. We write $\mathbf{W}_I^{\rm proj}$ (resp. $\mathbf{W}_I^{\rm inj}$) for the class of projective (resp. injective) $I$-weak equivalences. 
    \end{df}
    \begin{df} A morphism $p\colon \kX\to \kY$ is called a \textit{projective $I$-fibration} (resp.~\textit{an injective $I$-fibration}) if it has the right lifting property with respect to all maps that are both projective cofibrations (resp. monomorphisms) and projective (resp. injective) $I$-weak equivalences.
    \end{df}
Specializing Hirschhorn's Theorem \cite[Theorem 4.1.1]{Hirsch} to the case $\cM = \sPsh(T)_{\rm inj}$ or $\sPsh(T)_{\rm proj}$ with $S$ given by $\mathbf{W}_I^{\rm inj}$ or $\mathbf{W}_I^{\rm proj}$ respectively, produces the \textit{$I$-localized model structure} (injective or projective), denoted $L_I(\sPsh(T))$ with the relevant subscript. We will denote the homotopy category of $L_I(\sPsh(T))^{loc}_{\rm proj}$ (resp. $L_I(\sPsh(T))^{loc}_{\rm inj}$) by $\cH(T, I)^{loc}_{\rm proj}$ (resp. by $\cH(T,I)^{loc}_{\rm inj}$) or simply by $\cH(T,I)$ if it is clear wich model structure is considered. 
\medskip

The two classes of $I$-local projective and $I$-local injective weak equivalences agree, as proven in the following Proposition.
\begin{prop}\label{prop:proj-inj-local-agree}The identity functor $\id_{\sPsh(T)}\colon L_I(\sPsh(T))^{loc}_{\rm proj} \to L_I(\sPsh(T))^{loc}_{\rm inj} $ is the left adjoint of a Quillen equivalence.
    \begin{proof}    
   Recall first that the identity functor is a Quillen equivalence from the local projective to the local injective model structure before $I$-localization. Let $I^c\to I$ be the (functorially chosen) cofibrant replacement of $I$ in the projective local model structure on $\sPsh(T)$. Note in particular that $I^c\to I$ is a local weak equivalence between two injective-cofibrant objects. Write $\mathbf{W}^{\rm inj}_{I^c}$ for the class of $I^c$-injective weak equivalences, defined replacing $I$ with $I^c$ in Definitions \ref{def:I-local-proj-inj} and \ref{def:I-we}. We claim that $\mathbf{W}^{\rm inj}_{I^c} = \mathbf{W}^{\rm inj}_{I}$. Start with $f\colon A\to B \in \mathbf{W}^{\rm inj}_{I^c}$ and let $\kZ$ be an injective $I$-local object. Then $\kZ$ is globally fibrant (i.e.~fibrant for the injective local model structure) and thus (using that by Proposition \ref{prop:Day-makes-proj-loc-monoidal} $-\tensor \kX$ preserves local weak equivalences) the map $\cS(I\tensor \kX, \kZ) \to \cS(I^c\tensor \kX,\kZ)$ is a simplicial homotopy equivalence for every $\kX$ (note that since $\kZ$ is fibrant and $I^c\tensor \kX$ and $I\tensor \kX$ are both cofibrant, the simplicial mapping spaces are fibrant simplicial sets). We conclude that the map
    \[ \cS(I^c\tensor \kX,\kZ)\to  \cS(\kX, \kZ)\]
 is a weak equivalence. In particular, the object $\kZ$ is $I^c$-local and thus the map $f$ is an injective $I$-weak equivalence. 
We can reverse the argument, and start from $\kZ$ injective $I^c$-local to get that every injective $I$-weak equivalence is an injective $I^c$-weak equivalence. 

Note now that the classes of maps $\mathbf{W}^{\rm proj}_{I^c}$ and $\mathbf{W}^{\rm proj}_{I}$ do clearly coincide by definition. We are then reduced to show that $\mathbf{W}^{\rm proj}_{I^c} = \mathbf{W}^{\rm inj}_{I^c}$. This can be done using the same argument of \cite[Theorem 2.17]{DRO}. We recall the argument for completeness. Suppose that $f\colon A\to B$ is a projective $I^c$-weak equivalence and take $\kZ$ injective $I^c$-fibrant. Then $\kZ$ is globally fibrant, hence fibrant for the projective local model structure. The map of simplicial mapping spaces
\[ (f^c)^*\colon \cS(B^c, \kZ)\to \cS(A^c, \kZ) \]
is then a weak equivalence. Since the maps $A^c\to A$ and $B^c\to B$ are local weak equivalences, we conclude that $f\in  \mathbf{W}^{\rm inj}_{I^c}$. Conversely, start from $\kZ$ projective $I^c$-local and choose a local weak equivalence $\kZ\to \kZ'$ with $\kZ'$ globally fibrant (this exists by general principle, see \cite{JardineLocal} or \cite{JardineJPAA}). Then $\kZ'$ is easily seen to be injective $I^c$-fibrant. For every $f\colon A\to B$ injective $I^c$-weak equivalence, we then have a diagram of simplicial mapping spaces
\[\xymatrix{ \cS(B^c, \kZ) \ar[r] \ar[d]& \cS(A^c, \kZ) \ar[d] \\ 
\cS(B^c, \kZ') \ar[r] & \cS(A^c, \kZ') 
}\]
where the vertical arrows are weak equivalences and the bottom horizontal arrow is a weak equivalence since $f^c$  is an injective $I^c$-weak equivalence. Thus the top horizontal arrow is a weak equivalence as well, showing that $f$ is a projective $I^c$-weak equivalence. The fact that the identity functor gives a Quillen equivalence is now obvious. 
\end{proof}
    \end{prop}
\begin{rmk} For the category $\cM(k)$ of motivic spaces over $k$, a comparison between the two localized model structures (projective and injective) is given, for example, in \cite[Theorem 2.17]{DRO}: the identity functor $\id$ is the left adjoint of a Quillen equivalence between the motivic model structure (in the sense of \cite[2.12]{DRO}, built as localization of the projective model structure on simplicial presheaves) and the Goerss-Jardine model structure (built as localization of the injective model structure on simplicial presheaves).
    \end{rmk}
\subsubsection{}Let $f, g\colon \kX\to \kY$ be two morphisms of simplicial presheaves. As in \ref{subsec:Inverval-objects-monoidalcat}, an \textit{elementary $\Itensor$-homotopy from $f$ to $g$} is a morphism $H\colon \kX\tensor I\to \kY$ satisfying $H\circ \iota^I_0 = f$ and $H\circ \iota_1^I =g$. Two morphisms are called \textit{$\Itensor$-homotopic} if they can be connected by a sequence of elementary $\Itensor$-homotopies. A morphism $f\colon \kX\to \kY$ is called a \textit{strict $\Itensor$-homotopy equivalence} if there is a morphism $g\colon \kY\to \kX$ such that $f\circ g$ and $g\circ f$ are $\Itensor$-homotopic to the identity (of $\kY$ and of $\kX$ respectively). 
\begin{lem}[cfr. \cite{MV}, Lemma 3.6]\label{lem:I-homotopic-implies-Iwe} Any strict $\Itensor$-homotopy equivalence $f\colon \kX\to\kY$ is an $I$-weak equivalence for both the injective and the projective $I$-localized structure on simplicial presheaves.
    \begin{proof} We have to show that the compositions of $f$ with an $\Itensor$-homotopy inverse are equal to the corresponding identities in the homotopy categories  $\cH(T, I)_{\rm proj}$ and $\cH(T,I)_{\rm inj}$. But it's clear from the definition that two elementary $\Itensor$-homotopic maps are equal in the $I$-homotopy category (and this does not really depend on the choice of the injective/projective model structure).
        \end{proof}
    \end{lem}  
\subsection{Comparison of intervals}\label{sec:comparison-intervals} Let $\MM$ be the category of motivic spaces with modulus introduced in \ref{def:motivic-spaces}. 
Let $\cI = u_!v^*(\boxq)$ be the distinguished interval object of \ref{def:interval-I}. Let $\eta\colon \boxq\to \cI$ be the natural map constructed in \ref{subsubsec:ComparisonMapIntervals}.

Since $\boxq$ is a weak interval in $\MM$, we can talk about $\boxtensor$-homotopies between morphisms of motivic spaces with modulus.
\begin{prop}\label{prop:I-is-box-trivial} Let $\kX$ be any motivic space with modulus. Then the map $\id\tensor \iota_0^{\cI}\colon \kX = \kX\tensor \mathds{1} \to \kX\tensor \cI$ is a strict $\boxtensor$-homotopy equivalence.
    \begin{proof} Let $p = \id\tensor p_{\cI}\colon \kX\tensor \cI\to \kX\tensor \mathds{1} = \kX$ be the projection morphism. Since the composition $p\circ (\id\tensor \iota_0^{\cI})$ is the identity on $\kX$, it's enough to show that there exists a $\boxtensor$-homotopy between $(\id\tensor \iota_0^{\cI})\circ p$ and the identity on $\kX\tensor \cI$. Write $H$ for the map
        \[H= H_{\boxq}\colon (\kX \tensor \cI) \tensor \boxq \xrightarrow{\id_\kX\tensor\id_{\cI}\tensor \eta} \kX\tensor(\cI\tensor \cI) \xrightarrow{\id_{\kX}\tensor \mu} \kX\tensor \cI.\]
We need to check that $H\circ (\id\tensor \iota_0^{\boxq}) =(\id\tensor \iota_0^{\cI})\circ p$ and that $H\circ (\id\tensor \iota_1^{\boxq}) = \id$. By adjunction, the compositions $\eta \circ \iota_\epsilon^{\boxq}$ for $\epsilon = 0,1$ agree with the map $\iota_\epsilon^{\cI}$, so that $\id_{\kX\tensor \cI}\tensor \iota^{\boxq}_\epsilon = \id_{\kX\tensor \cI} \tensor \iota_\epsilon^{\cI}$. The required identities then follow from the interval structure on $\cI$.
        \end{proof}
    \end{prop}
Lemma \ref{lem:I-homotopic-implies-Iwe} and Proposition \ref{prop:I-is-box-trivial} together give the following
\begin{cor}\label{cor:map-X-to-XtensorI-is-box-we} For every motivic space $\kX$, the map $\kX\to \kX\tensor \cI$ is a $\boxq$-weak equivalence.
\end{cor}
\begin{cor}\label{cor:boxlocal-implies-Ilocal}Let $\kX$ be a $\boxq$-local object (for either the projective or the injective model structure on $\MM$). Then $\kX$ is $\cI$-local.
        \begin{proof}By definition, a   $\boxq$-local object $\kX$ satisfies the following condition: for every $\boxq$-weak equivalence $f\colon \kY\to \kZ$, the induced map $f^*$ on homotopy function complexes is a weak equivalence. In particular, the map
           \[\Map(\kY\tensor \cI, \kX)\to \Map(\kY, \kX)\]
           is a weak equivalence for every $\kY$, since $\kY\to \kY\tensor \cI$ is a $\boxq$-weak equivalence by Corollary \ref{cor:map-X-to-XtensorI-is-box-we}. But this is precisely the condition that a fibrant object $\kX$ has to satisfy for being $\cI$-local. 
    \end{proof}
\end{cor}
\begin{prop}Let $\kX$ be an $\A^1$-local motivic space with modulus. Then $\kX$ is $\cI$-local.
    \begin{proof}Let $\theta\colon \A^1\to \cI$ be the canonical map of \ref{subsubsec:canonical-maps-btw-intervals}, induced by adjunction by the identity morphism on $\A^1$ in $\Sm$. It's enough to show that the map $\id_\kX\tensor \iota_0^{\cI}\colon \kX\to \kX\tensor \cI$ is an $\A^1$ strict homotopy equivalence. An $\mathop{\mbox{$\A^1$-$\tensor$}}$-homotopy inverse is given by the projection $p\colon \kX\tensor \cI\to \kX$, and the homotopy between  $(\id_\kX\tensor \iota_0^{\cI})\circ p$ is the map
        \[H_{\A^1}\colon (\kX\tensor {\cI})\tensor \A^1 \xrightarrow{\id_{\kX\tensor {\cI}}\tensor \theta } \kX\tensor({\cI}\tensor {\cI})\xrightarrow{\id_\kX\tensor \mu} \kX\tensor {\cI}.\]
        The argument is then formally identical to the one given in Proposition \ref{prop:I-is-box-trivial} and Corollaries \ref{cor:map-X-to-XtensorI-is-box-we}  and \ref{cor:boxlocal-implies-Ilocal}.
        \end{proof}
    \end{prop}
We can also ask for the relation between the interval objects $v^*(\boxq)$ and $\A^1=(\A^1, \emptyset)$ in $\Mlog$. The functor $\omega_!$ is particularly well behaved, since it sends $\A^1$-local and $\boxq$-local objects in $\Mlog$ to $\A^1$-local objects in $\cM(k)$. In the other direction, we have the following
\begin{prop}Let $\kX$ be an $\A^1$-local motivic space in $\cM(k)$. Then $\lambda_*(\kX)$ is $v^*(\boxq)$-local in $\Mlog$.
    \end{prop}
    
\subsection{A singular functor}\label{sec:sing-functor} We specialize the results of the previous sections to $\MSmlog$, equipped with  the Nisnevich topology introduced in Section \ref{subsec:topology-Nis-Zar-modulus-cat}. This is the topology associated to a complete bounded regular cd-structure by Proposition \ref{prop:Nis-Zar-are-regu-comple-cd} and Proposition \ref{lem:Nis-Zar-are-bounded-cd}. In particular, we can apply Theorem \ref{thm:inj-loc-structure-is-bousfield}. 

\subsubsection{}\label{sec:review-properties-MM}Let $\MM_{\rm inj}^{loc}$ be the category of motivic spaces with modulus (over $k$) equipped with the local injective (for the Nisnevich topology) model structure and let ${\cI}$ be again the distinguished interval object of \ref{def:interval-I}. To simplify the notation, we write $\MM_{\rm inj}^{{\cI}-loc}$ for the ${\cI}$-localization of the local-injective model structure on $\MM$. By \cite[Theorem 4.1.1]{Hirsch}, $\MM_{\rm inj}^{\cI-loc}$ is a left proper cellular simplicial model category. Right properness does not follow formally.

For the category of motivic spaces $\cM(k)$ (without modulus), properness of the $\A^1$-local (injective) model structure is proved in \cite{MV}, Theorem 3.2 and, using a different technique, in \cite{JardineSpectra}, Theorem A.5. The proof of Morel and Voevodsky makes use of  the endofunctor  $\Sing_*$, that plays also an important role in the construction of a fibrant replacement functor.

Since we will work constantly in the category $\MM$, we will sometimes omit the locution ``with modulus'' for a motivic space $\kX$ (i.e.~for an object of $\MM$).

\subsubsection{}We introduce in this section an endofunctor $\SingI(-)$ on $\MM$ that plays in our theory the role of $\Sing_*$. Our results look formally like the corresponding statements in \cite[Section 3]{MV}, but the proofs are different. Another instance of this construction has been used in \cite[Appendix B]{BPO} in the context of motives for log schemes, where a similar problem (the lack of the multiplication map on $\overline{\square}$) arises. 

We start by noticing that the interval $\cI$ comes equipped with an extra diagonal map
\[\delta\colon \cI\to \cI\tensor \cI\]
induced by the diagonal $\delta = (\P^1, \infty)\to (\P^1\times \P^1, F_\infty^2)$ in  $\Smlog$. Thus, we can use the formulas of \ref{def:DeltaI} for constructing a cosimplicial object $\Delta_{\cI}^\bullet$ in $\MM$ whose $n$-th term is $\cI^{\tensor n}$. Similarly, we write $\Delta_{\A^1}^\bullet$ for the cosimplicial object deduced from $\A^1$ with the standard interval structure.
\begin{df}\label{def:functorSingI} Let $\kX$ be a motivic space with modulus. We write $\SingI(\kX)$ for the diagonal simplicial presheaf of the bi-simplicial presheaf
    \[ \Delta^{\rm op}\times \Delta^{\rm op}\to \Psh(\MSmlog), \quad ([n], [m])\mapsto [\Delta_{\cI}^m = \cI^{\tensor m}, \kX_n]. \]
   \end{df}
   For every $n\geq0$ there is canonical isomorphism $[\mathds{1}, \kX_n] = \kX_n$, giving by composition a natural transformation $s\colon \id \to \SingI$. Since the right adjoint $[-,-]=[-,-]^{\rm Day}$ to Day convolution preserves finite limits, for any $\kX$ the morphism $s_\kX\colon \kX\to \SingI(\kX)$ is a monomorphism and therefore a cofibration for the local-injective model structure.

\begin{prop}\label{prop:I-homotopy-implies-Sing-simplicial-homotopy}Let $f,g\colon \kX\rightrightarrows \kY$ be two morphisms of simplicial preshaves and let $H$ be an elementary $\cItensor$-homotopy between them. Then there exists an elementary simplicial homotopy between $\SingI(f)$ and $\SingI(g)$.
    \begin{proof}It is enough to show that there exists a simplicial homotopy
        \[\SingI(\kX)\times \Delta^1 \to \SingI(\kX\tensor \cI)\]
        between the natural maps $\SingI(\id_\kX\tensor \iota_0^{\cI})$ and $\SingI(\id_\kX\tensor \iota_1^{\cI})$. Thus, we have to construct a map of presheaves
         \[H_n\colon \SingI(\kX)_n\times \Delta^1_n (= \hom_\Delta([n],[1])) \to \SingI(\kX\tensor {\cI})_n = [{\cI}^{\tensor n}, \kX_n\tensor {\cI}]\]
         for every $n$, compatible with faces and degeneracies. Note here that ${\cI}$ is a discrete simplicial presheaf, so that $\kX_n\tensor {\cI} = (\kX\tensor {\cI})_n$ according to the definition of the extension of Day convolution to simplicial presheaves (see \eqref{eq:extension-Day-psh-to-spsh}). For every $M\in\MSmlog$, we have 
         \[ \SingI(\kX)_n(M) = \hom_{\Psh(\MSmlog)}(M\tensor {\cI}^{\tensor n}, \kX_n)\]
         so that a section over $M$ of $\SingI(\kX)_n$ is a map of presheaves $\alpha_M\colon M\tensor {\cI}^{\tensor n}\to \kX_n$.  Limits of presheaves are computed objectwise, and $\Delta^1_n$ is a constant presheaf, therefore a section over $M$ of $\SingI(\kX)_n\times \Delta^1_n$ is a pair $(\alpha_M, f)$ for $f\in \hom_\Delta([n], [1])$. Let $\Delta(f)\colon {\cI}^{\tensor n}\to {\cI}$ be the induced morphism given by the cosimplicial structure of $\Delta_{\cI}^\bullet$. Then we can consider the composition
         \begin{equation}\label{eq:def-homotopy-Sing} M\tensor {\cI}^{\tensor n}\xrightarrow{\id_M\tensor \delta_n} M\tensor {\cI}^{\tensor n} \tensor {\cI}^{\tensor n} \xrightarrow{\id_M\tensor {\cI}^{\tensor n} \tensor \Delta(f)} M\tensor {\cI}^{\tensor n}\tensor {\cI} \xrightarrow{\alpha_M \tensor \id_{\cI}} \kX_n \tensor {\cI}\end{equation}
         that defines a section over $M$ of $\SingI(\kX\tensor {\cI})_n = [{\cI}^{\tensor n}, \kX_n\tensor {\cI}]$. For $\phi\colon M'\to M$ there is a restriction map $\SingI(\kX)_n(M)\xrightarrow{\phi^*} \SingI(\kX)_n(M')$ that sends a section $\alpha_M$ to the composite $\alpha_M\circ (\phi\tensor \id)$. The assignment \eqref{eq:def-homotopy-Sing} is clearly compatible with $\phi^*$ and thus defines a morphism of presheaves of sets. It is easy to check that it is also compatible with the simplicial structure and that it gives indeed an homotopy between $\SingI(\id_\kX\tensor \iota_0^{\cI})$ and $\SingI(\id_\kX\tensor \iota_1^{\cI})$.
        \end{proof}
    \end{prop}  
\begin{cor}\label{cor:Sing-takes-i0-to-simpl-we}For any motivic space with modulus $\kX$, the morphism \[\SingI(\kX)\xrightarrow{\SingI(\id_{\kX}\tensor \iota_0^{\cI})} \SingI(\kX\tensor {\cI})\] is a simplicial homotopy equivalence (and the same holds for $\SingI(\id_{\kX}\tensor \iota_1^{\cI})$).
    \begin{proof} By Proposition \ref{prop:I-homotopy-implies-Sing-simplicial-homotopy}, it's enough to show that the map $\id_{\kX}\tensor \iota_0^{\cI}$ is an $\cItensor$-homotopy equivalence. A $\cItensor$-homotopy inverse is given by the projection map $p\colon \kX\tensor {\cI}\xrightarrow{\id\tensor p_{\cI}}\kX$. The composition $p\circ (\id_{\kX}\tensor \iota_0^{\cI})$ is clearly the identity, and an homotopy between $(\id_{\kX}\tensor \iota_0^{\cI})\circ p$ and the identity of $\kX\tensor {\cI}$ is  given by the multiplication map $\kX\tensor {\cI}\tensor {\cI}\xrightarrow{\id\tensor \mu} \kX\tensor {\cI}$.  The proof for $\SingI(\id_{\kX}\tensor \iota_1^{\cI})$ is similar, and it is left to the reader.
        \end{proof}
    \end{cor} 
\begin{lem}\label{lem:Xto-inhom-IX-we}The map $p^*\colon \kX = [\mathds{1}, \kX] \to [\cI, \kX]$ induced by the projection $p_{\cI}\colon \cI\to \mathds{1}$ is a strict $\cItensor$-homotopy equivalence.
    \begin{proof} A homotopy inverse is given by $\iota_0^*\colon [\cI, \kX]\to [\mathds{1},\kX]$ induced by $\iota_0^\cI$. Since the composition $\iota_0^*\circ p^*$ is the identity on $\kX$ we just need to show (as in the previous corollary) that there is a $\cItensor$-homotopy between the composition in the other direction $p^*\circ \iota_0^*$ and the identity morphism. By adjunction, we have a canonical isomorphism $[\cI\tensor \cI, \kX] \cong [\cI, [\cI, \kX]]$ that gives
        \[\hom_{\MM}([\cI, \kX], [\cI\tensor \cI, \kX]) \cong \hom_{\MM}([\cI, \kX]\tensor \cI, [\cI, \kX]),\]
        so that getting a homotopy $[\cI, \kX]\tensor \cI\to [\cI, \kX]$ it's equivalent to specifying a map $[\cI, \kX]\to[\cI\tensor \cI, \kX]$ satisfying the required identities. The multiplication map $\cI\tensor \cI\to \cI$ does the job.
        \end{proof}
    \end{lem}
\begin{prop}\label{prop:X-to-singX-is-we}For any motivic space with modulus $\kX$, the morphism $s_{\kX}\colon \kX\to \SingI(\kX)$ is an $I$-weak equivalence.
    \begin{proof} This is now formally identical to \cite[Corollary 3.8]{MV}, using our Lemma \ref{lem:Xto-inhom-IX-we}.
        \end{proof}
    \end{prop}
    The functor $\SingI$ has formally a left adjoint that is constructed, as usual, by left Kan extension. More precisely, we recall the following definition that is valid in every category of simplicial presheaves on a small site.
    \begin{df} Let $D^\bullet$ be a cosimplicial object in $\sPsh(T)$. We denote by $|-|_{D^\bullet}$ the left Kan extension ${\rm Lan}_{\Delta^\bullet}(D^\bullet)(-)$ of $D^\bullet$ along the functor $\Delta\to \sPsh(T)$ that sends $[n]$ to the constant simplicial presheaf $\Delta^n$.
        \end{df}
There is a canonical isomorphism $|\Delta^n|_{D^\bullet} = D^n = D([n])$. We can write an explicit description as follows. Let $\kX$ be a simplicial presheaf and identify every $\kX_n$ with a simplicial presheaf of dimension zero. Then $|\kX|_{D^\bullet}$ is the co-equalizer (in $\sPsh(T)$) of the following diagram
    \begin{equation}\label{eq:realization-is-coequ}
    \xymatrix{
  \coprod_{\phi\in \hom_\Delta([m],[n])} \kX_n\times D^m   \ar@<.5ex>[r]^-{f} \ar@<-.5ex>[r]_-{g} & \coprod_{n}\kX_n\times D^n
    }
    \end{equation}
    where the $n$-th term of the map $f\colon  \kX_n\times D^m  \to  \kX_n\times D^n$ is induced by the cosimplicial structure of $D^\bullet$  and the $m$-th term of the map $g\colon  \kX_n\times D^m  \to  \kX_m\times D^m$ is induced by the simplicial structure of $\kX$.
  \begin{rmk}   We specialize to the case $D= \Delta^\bullet \tensor \Delta_{\cI}^\bullet$ in $\MM$. By definition, the functor $\SingI(-)$ satisfies, for any $\kY\in \MM$,
    \begin{align*}\SingI(\kY)_n & = \hom_{\MM}(\Delta^n, \SingI(\kY)) = \hom_{\MM}(\Delta^n,[{\cI}^{\tensor n} , \kY] ) \\ &= \hom_{\MM}(\Delta^n \tensor {\cI}^{\tensor n}, \kY)  = \hom_{\MM}(|\Delta^n|_{\Delta^\bullet\tensor \Delta_{\cI}^\bullet}, \kY)\end{align*}
        so that it's clear by construction that $|-|_{\Delta^\bullet\tensor \Delta_{\cI}^\bullet}$ is its left adjoint.
            \end{rmk}
We summarize the properties of the functor $\SingI(-)$ proved so far in the following Theorem.
\begin{thm}\label{thm:properties-SingI}The endofunctor $\SingI(-)$ of $\MM$ commutes with limits and takes the morphism $\id\tensor \iota_0^{\cI}\colon \kX\to \kX\tensor \cI$ to a simplicial weak equivalence for every motivic space $\kX$. Moreover, the canonical natural transformation $\id\to \SingI(-)$ is both a monomorphism and an injective $\cI$-weak equivalence.
    \begin{proof}  $\SingI(-)$ commutes with limits since it is by construction a right adjoint. The other statements are precisely the content of Corollary \ref{cor:Sing-takes-i0-to-simpl-we} and Proposition \ref{prop:X-to-singX-is-we}.
        \end{proof}
    \end{thm}
By the same argument discussed after Definition \ref{def:I-local-proj-inj}, the functor $\SingI(-)$ takes also the morphisms $\id\tensor p_{\cI} \colon \kX\tensor \cI\to \kX$ to a simplicial weak equivalence for every motivic space $\kX$. Applying the argument again, the same holds for the morphisms $\id\tensor \iota_1^{\cI}$. 

\medskip

To continue or construction, we need to further study the properties of the adjoint    $|-|_{\Delta^\bullet\tensor \Delta_{\cI}^\bullet}$ to $\SingI(-)$.         
    \begin{rmk} For $D^\bullet = \Delta^\bullet$, we have ${\rm Lan}_{\Delta^\bullet}(D^\bullet)(-) = {\rm Lan}_{\Delta^\bullet}(\Delta^\bullet)(-) \cong \id$, that is, $|-|_{\Delta^\bullet}$ is the (pointwise) left Kan extension of the functor $\Delta^\bullet$ along itself and is isomorphic to the identity functor. 
        \end{rmk}
\begin{lem}[see \cite{MV}, Lemmas 3.9 and 3.10]\label{lem:realization-preserves-mono}Let $D^\bullet$ be a cosimplicial object in $\MM$. Then the following conditions are equivalent
    \begin{romanlist}
        \item the morphism $D^0\amalg D^0\to D^1$ induced by the cofaces is a monomorphism,
        \item the functor $|-|_{D^\bullet}$ preserves monomorphism.
        \end{romanlist}
\end{lem}
Lemma \ref{lem:realization-preserves-mono} clearly applies to $D^\bullet = \Delta^\bullet\tensor {\cI}^\bullet$. In this case we have $D^0 = \Delta^0\tensor {\cI}^{\tensor 0} = \Delta^0 \tensor \mathds{1} = \pt$ and the maps $D^0\rightrightarrows D^1 = \Delta^1\tensor {\cI}$ induced by the cofaces are two distinct rational points $\pt\rightrightarrows \Delta^1\tensor {\cI}$. 
\begin{lem}\label{lem:tensor-and-prod}Let $\kX$ and $\kY$ be two motivic spaces. Then the natural map $\phi\colon (\kX\times \kY)\tensor {\cI}\to \kX\times (\kY \tensor {\cI})$ is an ${\cI}$-homotopy equivalence, and hence an ${\cI}$-weak equivalence.
    \begin{proof} The existence of the map is guaranteed by universal property given the existence of morphisms $(\kX\times \kY)\tensor {\cI} \xrightarrow{p_1\tensor \id} \kX\tensor {\cI}\xrightarrow{\id \tensor p_I} \kX$ and $(\kX\times \kY)\tensor {\cI} \xrightarrow{p_2\tensor \id}\kY\tensor {\cI}$ for $p_1, p_2$ the projections $\kX\times \kY\to \kX$ and $\kX\times \kY\to \kY$.  Let $\nu$ be the morphism
        \[ \kX\times (\kY \tensor {\cI}) \xrightarrow{\id \times(\id \tensor p_{\cI})} \kX\times \kY \xrightarrow{\id\tensor \iota_0^{\cI}} (\kX\times kY)\tensor {\cI}.\]
We claim  that $\nu \circ \phi$ is $\cItensor$-homotopic to the identity of $(\kX\times \kY)\tensor {\cI}$ and that $\phi\circ \nu$ is $\cItensor$-homotopic to the identity of $\kX\times(\kY\tensor {\cI})$. For the first composition, a homotopy is given by the multiplication map 
\[(\kX\times \kY) \tensor {\cI}\tensor {\cI}\xrightarrow{\id\tensor \mu} (\kX\times \kY) \tensor {\cI}\]
and the required identities follow from the fact that the morphism $(\kX\times \kY)\tensor {\cI}\to \kX\times \kY$ factors through $\phi$. For the second composition, we consider the morphism $H$ defined as
\[ (\kX\times (\kY \tensor {\cI})) \tensor {\cI} \to \kX \times ((\kY\tensor {\cI})\tensor {\cI}) = \kX \times (\kY\tensor ({\cI}\tensor {\cI})) \xrightarrow{\id\times (\id\tensor \mu)} \kX\times (\kY \tensor {\cI})\]
where the first map is given again by universal property. It is easy to check that this is indeed the required homotopy.
        \end{proof}
\end{lem}
\begin{rmk}Let $A$ be the class of morphisms $\{\kX\xrightarrow{\id\tensor \iota_0^{\cI}} \kX\tensor {\cI}\}_{\kX\in \MM}$. Following \cite[Definition 2.1]{MV}, we say that a local injective fibrant motivic space $\kY$ is $A$-local if for any $\kZ$ in $\MM$ and any $\kX\xrightarrow{\id\tensor \iota_0^{\cI}} \kX\tensor {\cI}$ in $A$, the map
    \[\Map(\kZ\times(\kX\tensor {\cI}), \kY)\to \Map(\kZ\times \kX, \kY)\]
    is a weak equivalence of simplicial sets. By Lemma \ref{lem:tensor-and-prod}, an object $\kZ$ is $A$-local if and only if it is (injective) ${\cI}$-local according to Definition \ref{def:I-local-proj-inj}.
    \end{rmk}
\begin{lem}\label{lem:realization-preserves-I-we} For any motivic space with modulus $\kX$, the morphism $|\kX|_{\Delta^\bullet\tensor \Delta_{\cI}^\bullet} \to |\kX|_{\Delta^\bullet} \cong \kX$ induced by the projection $p_{\cI}\colon \Delta^\bullet\tensor \Delta_{\cI}^\bullet \to \Delta^\bullet$ is an ${\cI}$-weak equivalence.
    \begin{proof} The functor $|-|_{\Delta^\bullet\tensor \Delta_{\cI}^\bullet}$ commutes with colimits and we can do induction on the skeleton to the reduce to the case $\kX = \kY \times \Delta^r$ for some  simplicial presheaf $\kY$ of simplicial dimension $0$ and some $r\geq 0$. The co-equalizer \eqref{eq:realization-is-coequ} then takes the form
        \[\xymatrix{
  \coprod_{\phi\in \hom_\Delta([m],[n])} \kY \times \Delta^r[n]\times (\Delta^m\tensor {\cI}^{\tensor m})   \ar@<.5ex>[r]^-{f} \ar@<-.5ex>[r]_-{g} & \coprod_{n}\kY \times \Delta^r[n]\times (\Delta^n\tensor {\cI}^{\tensor n}). 
    }\]
    The term $\kY$ is constant and the functor $\kY\times (-)$  commutes with arbitrary colimits in $\MM$. So we can identify $|\kY\times \Delta^r|_{\Delta^\bullet\tensor {\cI}^\bullet}$ with $\kY\times (\Delta^r\tensor {\cI}^{\tensor r})$. The morphism $|\kY \times \Delta^r|_{\Delta^\bullet\tensor \Delta_{\cI}^\bullet}\to \kY \times \Delta^r$ is then identified with the projection $\kY \times (\Delta^r \tensor {\cI}^{\tensor r}) \to \kY\times \Delta^r$. By Lemma \ref{lem:tensor-and-prod} (applied several times), the product $\kY \times (\Delta^r \tensor {\cI}^{\tensor r})$ is ${\cI}$-weakly equivalent to the product $(\kY \times \Delta^r) \tensor {\cI}^{\tensor r}$, and the projection $p_{\cI}^r\colon (\kY \times \Delta^r) \tensor {\cI}^{\tensor r}\to \kY \times \Delta^r$ factors through the canonical map $(\kY \times \Delta^r) \tensor {\cI}^{\tensor r}\to \kY \times (\Delta^r \tensor {\cI}^{\tensor r})$. An easy induction shows that $p_{\cI}^r$ is an ${\cI}$-weak equivalence, and the statement follows.
        \end{proof}
    \end{lem}
\begin{prop}\label{prop:SingI-preserves-fibrations}The functor $\SingI(-)$ preserves injective ${\cI}$-fibrations.
    \begin{proof}Equivalently (using the same strategy of \cite{MV}, Corollary 3.13), we show that its left adjoint $|-|_{\Delta^\bullet\tensor \Delta_{\cI}^\bullet}$ preserves injective cofibrations (i.e. monomorphisms of presheaves) and injective ${\cI}$-weak equivalences. The first property is provided by Lemma \ref{lem:realization-preserves-mono}, while the second property is provided by Lemma \ref{lem:realization-preserves-I-we} together with the 2 out of 3 property of injective ${\cI}$-weak equivalences.
        \end{proof}
    \end{prop}
\subsection{Properness}\label{sec:properness} We can now prove right properness of the injective ${\cI}$-local model structure $\MM_{\rm inj}^{loc}$. The proof is a combination of arguments due to Jardine in \cite[Appendix A]{JardineSpectra} and Morel-Voevodsky \cite[pp.~77--82]{MV}.
\begin{lem}\label{lem:can-apply-properness} Let $\kX,\kY, \mathscr{E}$ be three motivic spaces with modulus. Let $p\colon \mathscr{E}\to \kX\times (\kY\tensor {\cI})$ be an injective ${\cI}$-fibration. Then in the cartesian square
    \begin{equation}\label{eq:lem-can-apply-properness}\xymatrixcolsep{5pc}\xymatrix{ \mathscr{E}\times_{\kX\times (\kY\tensor {\cI})} (\kX\times \kY) \ar^-f[r] \ar[d]& \mathscr{E}\ar[d]^p \\
    \kX\times \kY \ar[r]^-{\id_{\kX}\times (\id_{\kY}\tensor \iota_0^{\cI})} & \kX\times(\kY \tensor {\cI})
    }\end{equation}
    the morphism $f$ is an injective ${\cI}$-weak equivalence.
    \begin{proof}This statement is proved by using a combination of the established properties of the functor $\SingI(-)$. Applying $\SingI(-)$ to the top arrow of \eqref{eq:lem-can-apply-properness} gives the following diagram
        \[\xymatrixcolsep{4pc} \xymatrix{
        \SingI(\mathscr{E}\times_{\kX\times (\kY\tensor {\cI})} (\kX\times \kY)) \ar[r]^-{\SingI(f)} & \SingI(\mathscr{E}) \\
        \mathscr{E}\times_{\kX\times (\kY\tensor {\cI})} (\kX\times \kY) \ar[u] \ar[r]^-f & \mathscr{E}\ar[u]
        }
        \]
        whose vertical maps are injective ${\cI}$-weak equivalences by Theorem \ref{thm:properties-SingI}. We will show that $\SingI(f)$ is a simplicial weak equivalence. Note that again by Theorem \ref{thm:properties-SingI}, $\SingI(-)$ commutes with limits, so that applying it to \eqref{eq:lem-can-apply-properness} gives another cartesian square. Using Lemma \ref{lem:tensor-and-prod} and Theorem \ref{thm:properties-SingI}, we see that $\SingI(\kX\times \kY)\to \SingI(\kX\times (\kX\tensor {\cI}))$ is a simplicial weak equivalence. By Proposition \ref{prop:SingI-preserves-fibrations}, the morphism $\Sing^I(\mathscr{E})\to \SingI(\kX\times(\kY\tensor {\cI}))$ is a simplicial fibration. Since the injective model structure on simplicial presheaves is proper, we conclude that $\SingI(f)$ is a weak equivalence as required.
        \end{proof}
\end{lem}
\begin{lem}\label{lem:replace-Jardine-A1}Suppose we are given morphisms of motivic spaces $\kX\xrightarrow{\id\tensor \iota_0^{\cI}} \kX\tensor {\cI} \xrightarrow{g} \kY$ and an injective ${\cI}$-fibration $p\colon \mathscr{E}\to \kY$. Then the induced map
    \[ \mathscr{E}\times_{\kY} \kX\to \mathscr{E}\times_{\kX} (\kX\tensor {\cI})\]
    is an injective ${\cI}$-weak equivalence.
    \begin{proof}The class of fibrations in $\MM_{\rm inj}^{loc}$ is closed under pull-backs (as in any closed model category), so that the morphism  $\mathscr{E}\times_{\kX} (\kX\tensor {\cI})\to \kX\tensor {\cI}$ is an injective ${\cI}$-fibration. The statement then follows from Lemma \ref{lem:can-apply-properness}.
        \end{proof}
    \end{lem}
\begin{lem}\label{lem:replace-corA4-Jardine}Suppose we are given a morphism of motivic spaces $\kX\xrightarrow{g} \kY$ and an  injective ${\cI}$-fibration $p\colon \mathscr{E}\to \kY$ with $\kY$ fibrant for the injective ${\cI}$-local model structure. Suppose moreover that $g$ is an injective ${\cI}$-weak equivalence. Then the top horizontal arrow in the pull-back square
    \[\xymatrix{ \mathscr{E}\times_{\kY} \kX \ar[r]\ar[d]  & \mathscr{E}\ar[d]^-p\\
    \kX\ar[r]^g & \kY
    }\]
    is an injective ${\cI}$-weak equivalence.
    \begin{proof} Choose a factorization of $g$ as $\kX\xrightarrow{j} \mathscr{W} \xrightarrow{q}$ with $q$ an injective ${\cI}$-fibration and $j$ an elementary ${\cI}$-cofibration (see \cite[p.~75]{MV} where the class $B_1$ of elementary $A$-cofibrations is introduced  or \cite[p.~537]{JardineSpectra}). Since $\kY$ is fibrant, $\mathscr{W}$ is ${\cI}$-fibrant as well. Since $g$ and $j$ are injective ${\cI}$-weak equivalences, $q$ is an injective ${\cI}$-weak equivalence between ${\cI}$-fibrant objects, and it is therefore a local weak equivalence. Since $p$ is a global fibration and the injective local model structure is proper, $q$ pulls back along $p$ to a local weak equivalence (and thus to an ${\cI}$-weak equivalence) $\mathscr{W}\times_{\kY} \mathscr{E} \to \mathscr{E}$. We are then left to show that the natural map $\mathscr{E}\times_{\kY} \kX \to \mathscr{W}\times_{\kY} \mathscr{E}$ is an $I$-weak equivalence. This follows from \cite[Lemma A.3]{JardineSpectra}, using Lemma \ref{lem:replace-Jardine-A1} instead of  \cite[Lemma A.1]{JardineSpectra}.
        \end{proof}
    \end{lem}
Write $\phi_{\kX}\colon \kX\to \mathcal{L}\kX$ for a functorial ${\cI}$-injective fibrant model of a motivic space $\kX$. The morphism $\phi_{\kX}$ is by construction a cofibration and an injective ${\cI}$-weak equivalence, while the map $\mathcal{L}\kX\to \pt$ is an injective ${\cI}$-fibration. Note that the existence of $\mathcal{L}\kX$ is guaranteed by the fact that $\MM_{\rm inj}^{loc}$ is cofibrantly generated.
\begin{thm}\label{thm:properness}Let $\kX, \kY$ and $\mathscr{E}$ be motivic spaces and suppose we are given a diagram
    \[\xymatrix{ \mathscr{E}\times_{\kY} \kX \ar[r]^-\gamma\ar[d]  & \mathscr{E}\ar[d]^-p\\
    \kX\ar[r]^g & \kY
    }\]
    where $g$ is an injective ${\cI}$-weak equivalence and $p$ is an injective ${\cI}$-fibration. Then the induced morphism $ \gamma\colon \mathscr{E}\times_{\kY} \kX \to \mathscr{E}$ is also an injective ${\cI}$-weak equivalence.
    \begin{proof}We can construct a commutative square of the form
        \[ \xymatrix{ \mathscr{E}\ar[r]^j \ar[d]^-p &  \mathscr{E}'\ar[d]^-{p'} \\
        \kY \ar[r]^{\phi_{\kY}} &\cL \kY
        }\]
such that the upper horizontal arrow $j$ is a cofibration and an ${\cI}$-weak equivalence and $p'$ is an injective ${\cI}$-fibration. By Lemma \ref{lem:replace-corA4-Jardine}, the induced map $\phi'\colon \mathscr{E}'\times_{\cL \kY} \kY \to \mathscr{E}'$ is an injective ${\cI}$-weak equivalence between ${\cI}$-fibrant objects over $\kY$, and therefore the morphism $\theta\colon \mathscr{E}\to \mathscr{E}'\times_{\cL \kY} \kY$ is also an injective ${\cI}$-weak equivalence. Since both objects are  fibrant over $\kY$,  $\theta$ is also a simplicial homotopy equivalence, and so also a local weak equivalence (cfr.~\cite[Lemma 2.30]{MV}). In the cartesian square
\[ \xymatrix{ \kX\times_{\kY} \mathscr{E} \ar[r]^-{\gamma} \ar[d]^{\theta'} & \mathscr{E}\ar[d]^\theta \\
\kX\times_{\cL \kY} \mathscr{E}' \ar[r]^{\gamma'} & \kY\times_{\cL\kY} \mathscr{E}'
}
\]
the morphism $\gamma'$ is an injective ${\cI}$-weak equivalence, since $\phi'\circ \gamma'$ is an injective ${\cI}$-weak equivalence by Lemma \ref{lem:replace-corA4-Jardine} (and we have already noticed that $\phi'$ is an injective ${\cI}$-weak equivalence). We are then left to show that $\theta'$ is an injective ${\cI}$-weak equivalence to conclude. But this map is in fact a local weak equivalence, since the morphism $\theta$ is one and the injective local model structure on simplicial presheaves is proper.
\end{proof} 
\end{thm}
\begin{rmk}Theorem \ref{thm:properness} is  a reformulation in our context of \cite[Theorem A.5]{JardineSpectra}. This is a consequence of a more general statement, as explained in \cite[7.3]{JardineLocal}. Jardine proves in \textit{loc.~cit.~}Lemma 7.25, right properness for the $f$-local theory where $f\colon A\to B$ is a cofibration of simplicial presheaves, satisfying certain conditions, analogue to the property established by Lemmas \ref{lem:replace-Jardine-A1} and \ref{lem:replace-corA4-Jardine}. However, the category $\MM_{\rm inj}^{{\cI}-loc}$ is not obtained as $f$-local theory (in the sense of Jardine), since we are inverting $\kX\to \kX\tensor {\cI}$ and not $\kX\to \kX \times {\cI}$. This explains why, for example, the proof of \cite[Lemma A.1]{JardineSpectra} does not go through in our context and we need to follow more closely \cite{MV}, using the singular functor $\SingI(-)$. From this point of view, Theorem \ref{thm:properness} is new and does not follow from Jardine's work.
\end{rmk}
\begin{rmk}\label{rmk:proj-is-right-proper}Having Theorem \ref{thm:properness} at hand, it is possible to use the strategy of \cite[Lemma 3.4]{Blander} to show that the ${\cI}$-local projective model structure on $\MM$ is also right proper. The proof in \textit{loc.cit.}~uses in an essential way the right properness of the injective structure. We leave the details to the interested reader.
    \end{rmk}
\subsection{Nisnevich B.G.~property}\label{sec:Nis-BG} Let $\MM^{{\cI}-loc}_{\rm inj}$ be again the category of motivic spaces with modulus equipped with the ${\cI}$-local injective model structure. To shorten the notation, we will refer to it as the \textit{(injective) motivic model structure}. One can define in a similar way a projective variant. We set the following notation.

\begin{df}The class of ${\cI}$-local-weak equivalences will be called the class of \textit{motivic weak equivalences}. We say that a map $f\colon \kX\to \kY$ in $\MM$ is  \textit{motivic fibration}  if it is an injective ${\cI}$-fibrations. An object $\kX$ is \textit{motivic fibrant} if the structure morphism is a motivic fibration.\end{df}
Note that by Proposition \ref{prop:proj-inj-local-agree} we don't need to specify if we consider the injective or the projective model structure in the definition of motivic weak equivalences.

 The \textit{unstable unpointed motivic homotopy category with modulus} $\MH$ over $k$ is the homotopy category associated to the model structure $\MM^{{\cI}-loc}_{\rm inj}$. 

For $M\in \MSmlog$ and $\kX\in \MM$, we write $\kX(M\tensor {\cI})$ for the simplicial set $\Map(M\tensor {\cI}, \kX)$. The morphism $\id\tensor \iota_0^{\cI}$ induces then a morphism a simplicial sets $\kX(M\tensor {\cI}) \to \kX(M)$.
\begin{lem}\label{lem:Bousfield-loc-description-fibrant-motivic}A motivic space $\kX \in \MM$ is motivic fibrant if and only if 
    \begin{enumerate}\item $\kX$ is globally fibrant (see \cite[5.1]{JardineLocal}), and
        \item for every $M\in \MSmlog$, the map  $\kX(M\tensor {\cI}) \to \kX(M)$ is a weak equivalence.
        \end{enumerate} 
    \end{lem}
    
\begin{df}A motivic space with modulus is called ${\cI}$-invariant if for every $M\in \MSmlog$, the map  $\kX(M\tensor {\cI}) \to \kX(M)$ is a weak equivalence.
    \end{df}
    \begin{df}\label{def:Nisn-excisive-prop}A motivic space with modulus $\kX$ is said to be \textit{Nisnevich excisive} (or to have the \textit{B.G.~property with respect to Nisnevich squares}) if $\kX(\emptyset)$ is contractible and the square
\[\xymatrix{ \kX(\ol{M}; \partial M, D_M) \ar[r]^{j^*} \ar[d]^{p^*} & \kX(\ol{U}; \partial U, D_U)\ar[d] \\
\kX({\ol{Y}}; \partial Y, D_Y) \ar[r] &  \kX(\ol{U} \times_{\ol{M}} \ol{Y}, \partial Y \cap p^{-1}({\ol{U}}), D_{\ol{U} \times_{\ol{M}} \ol{Y}})
}
\]
is homotopy cartesian for every elementary Nisnevich square of the form \eqref{eq:Nis-square-MSmlog}.
        \end{df}
The following Proposition is the analogue of \cite[Proposition 1.16]{MV} in our setting.
\begin{prop}[B.G.~property for motivic spaces]\label{prop:BG}Let $\kX$ be a motivic space with modulus. Then the following are equivalent
    \begin{enumerate}
        \item\label{Nis-ex-and-I-invariant} $\kX$ is Nisnevich excisive  and ${\cI}$-invariant
        \item Any motivic fibrant replacement $\kX\to \cL \kX$ (i.e.~a fibrant replacement for the ${\cI}$-injective model structure) is a sectionwise weak equivalence. 
        \end{enumerate}
    \begin{proof}We follow the proof of  \cite[Proposition 1.16]{MV}. The fact that the second condition implies the first one follows from the explicit description of fibrant objects in a Bousfield localization. Indeed, it's clear that if  $\kX\to \cL \kX$ is a sectionwise weak-equivalence, then $\kX$ is ${\cI}$-invariant, because if the map $\kX\to \cL \kX$ is a sectionwise weak equivalence, for every motivic space $\kZ$ the induced map $\Map(\kZ, \kX)\to \Map(\kZ, \cL \kX)$ is weak equivalence of simplicial sets, and this applies in particular to $\kZ = M\tensor {\cI}$. To show the B.G.~property, we combine Lemma \ref{lem:Bousfield-loc-description-fibrant-motivic} with Theorem \ref{thm:inj-loc-structure-is-bousfield}, noting that the B.G.~property is invariant with respect to weak equivalences of presheaves i.e.,~sectionwise weak equivalences (see \cite[Remark 3.1.14]{MV}).
        
        Conversely, suppose that $\kX$ satisfies condition (\ref{Nis-ex-and-I-invariant}). Let $\phi_{\kX}\colon \kX\to \cL \kX$ be a motivic fibrant replacement. Using \cite[Proposition 3.8]{Vcdtop} we are left to show that $\phi_{\kX}$ is a local injective fibrant replacement to deduce that it is a sectionwise weak equivalence.  Write a factorization of $\phi_{\kX}$ as \[\kX\xrightarrow{\psi} \kY \xrightarrow{\phi'} \cL \kX\]in the local injective model structure with $\psi$ a trivial local cofibration (so a monomorphism that is also a local weak equivalence) and $\phi'$ a local injective fibration. Note that $\psi$ is a motivic weak equivalence, so that also $\phi'$ is a motivic weak equivalence. Since $\phi'$ is a local injective fibration and $\cL \kX$ is local injective fibrant, $\kY$ is also local injective fibrant. Hence $\psi$ is a trivial local cofibration with target a local injective fibrant object, and therefore it is a fibrant replacement for $\kX$ in the local injective model structure. In particular, $\kY$ is local projective fibrant, hence by Proposition \ref{prop:BlanderProp} it has the B.G.~property with respect to Nisnevich squares (it is flasque in the sense of \cite[Definition 3.3]{Vcdtop}). By \cite[Lemma 3.5]{Vcdtop} (see also \cite[Lemma 1.18]{MV}), the map $\psi$ is then a sectionwise weak equivalence and so by Lemma \ref{lem:Bousfield-loc-description-fibrant-motivic} we conclude that $\kY$ is motivic fibrant. Since $\phi'$ is now an injective ${\cI}$-weak equivalence between injective ${\cI}$-fibrant objects, we conclude that $\phi'$ is an objectwise simplicial weak equivalence by \cite[Theorem 3.2.13.(1)]{Hirsch}.
        \end{proof}
    \end{prop}
    \begin{cor}  Any morphism $f\colon \kX\to\kY$ between motivic spaces with modulus satisfying the equivalent conditions of Proposition \ref{prop:BG} is a motivic weak equivalence if and only if it is a sectionwise weak equivalence.
        \begin{proof}Taking motivic fibrant replacements of $\kX$ and of $\kY$ we get a diagram
        \[\xymatrix{ \kX\ar[r]^f \ar[d]_{\phi_{\kX}} & \kY \ar[d]^{\phi_{\kY}}\\
        \cL \kX \ar[r]^{f'} & \cL \kY
         }
        \]
        with vertical arrows sectionwise weak equivalences. The induced morphism $f'$ is a motivic weak equivalence if (and only if) $f$ is a motivic weak equivalence. The local Whitehead Lemma \cite[Theorem 3.2.13]{Hirsch} implies that $f'$ is  an objectwise weak equivalence in this case, and thus so is $f$.  
            \end{proof}
        \end{cor}
        
        \subsection{Pointed variants}\label{sec:pointed-spaces} We write $\MMb$ for the category of \textit{pointed motivic spaces with modulus}. Objects are pairs, $(\kX, x)$, where $\kX$ is a motivic space with modulus and $x\colon \pt\to \kX$ is a fixed basepoints. Equivalently, a pointed motivic space with modulus is a contravariant functor 
        \[\kX\colon \MSmlog\to \cS_\bullet\]
        from $\MSmlog$ to the category $\cS_\bullet$ of pointed simplicial sets. We have a canonical adjunction
        \[ (-)_+\colon \MM\leftrightarrows\MMb \]
        where $(-)_+$ is the ``add base point functor'', left adjoint to the forgetful functor. The category $\MMb$ inherits by \cite[7.6.5]{Hirsch} two natural model structures, induced respectively by the $I$-local injective and projective model structure on $\MM$.
        \begin{thm}\label{thm:pointed-version-model-cats} The injective ${\cI}$-local structure on $\MMb$ is proper, cellular and simplicial, while the projective ${\cI}$-local structure on $\MMb$ is left proper, cellular and simplicial (but see Remark \ref{rmk:proj-is-right-proper} on right properness). A morphism $f\colon (\kX, x) \to (\kY, y)$ is an injective (resp. projective) ${\cI}$-weak equivalence if and only if the underlying morphism of unpointed motivic spaces $\kX\to \kY$ is a weak equivalence in the ${\cI}$-localized injective (resp. projective) motivic model structure. The statement remains true if one replaces the projective ${\cI}$-localized model structure on $\MM$ with the ${\cI}^c$-localized projective structure on $\MM$.
            \end{thm}
  
            By Proposition \ref{prop:proj-inj-local-agree}, the identity functor is a (left) Quillen equivalence between the projective and the injective motivic model structure on $\MMb$. 
            We denote by $\MHp$ the \textit{unstable pointed motivic homotopy category with modulus} over $k$.
            
        
        \subsubsection{}There are two natural closed monoidal structures on $\MMb$, induced respectively by the cartesian product and by Day convolution on $\MM$ (see Remark \ref{rmk:pointed-Day}).
        \begin{romanlist}
            \item \textit{Smash product}. We define the smash product objectwise, $\kX\wedge \kY (M) = \kX(M)\wedge \kY(M)$ for $\kX$ and $\kY$ pointed motivic spaces. With this definition $\MMb$ is naturally enriched over $\cS_\bullet$. The simplicial function space is given degreewise by the pointed simplicial set $\cS_\bullet(\kX, \kY)$
            \[\cS_\bullet(\kX, \kY)_n = \hom_{\MMb}(\kX\wedge (\Delta[n])_+, \kY), \]
            and internal hom given by 
            \[\Hom_{\MMb}(\kX, \kY)(M) = \cS_\bullet(\kX, \kY)(\kX\wedge (h_M)_+, \kY).\]
            \item \textit{Pointed Day convolution}. For two pointed motivic spaces $\kX$ and $\kY$, we define their pointed convolution $\kX\Daytensor_\bullet \kY$ by means of the following push-out diagram of unpointed simplicial presheaves
            \[\xymatrix{ (\kX\Daytensor \pt)\coprod (\pt \Daytensor \kY) \ar[r] \ar[d] & \kX\Daytensor \kY \ar[d]\\
             \pt \ar[r] & \kX\Daytensor_\bullet \kY,
            }
            \]
            where the top horizontal arrow is a monomorphism. We have, in particular, $(\kX\tensor \kY)_+ = (\kX)_+ \Daytensor_\bullet (\kY)_+$ for unpointed motivic spaces $\kX$ and $\kY$. The unit for $\Daytensor_\bullet$ is $\mathds{1}_+$. This definition makes the add base point functor $(-)_+$ strict monoidal.  We denote by $[-,-]_\bullet$ the pointed version of the internal hom for Day convolution.
            \end{romanlist}
\begin{lem}\label{lem:smash-preserves-I-proj-we}Let $\MMb^{{\cI}^c-loc}_{\textit{proj}}$ denote the category of pointed motivic spaces equipped with the ${\cI}^c$-local projective model structure. Then taking the smash product with any motivic space $\kX$ preserves projective ${\cI}^{c}$-weak equivalences.
    \begin{proof}The argument given in \cite[Lemma 2.18 - Lemma 2.20]{DRO} works almost \textit{verbatim} in our setting, so we just sketch the proof. 
        First, we prove that smashing with any cofibrant space $\kX$ preserves projective ${\cI}^{c}$-weak equivalences.
        Given a pointed motivic space $\kZ$ that is ${\cI}^c$-fibrant and a pointed motivic space $\kX$ that is cofibrant, the internal hom $\Hom_{\MMb}(\kX, \kZ)$ is clearly objectwise fibrant since the category of pointed simplicial presheaves on any small site $T$ equipped with the projective structure is a monoidal model category for the smash product. 
        In particular the internal hom $\Hom_{\MMb}(\kX, -)$ is a right Quillen functor. Let  now $\Lambda'$ be the set of maps 
        \[\Lambda' = \Sigma_P \cup \{ (h_M)_+ \Daytensor_\bullet {\cI}^c_+ \to (h_M)_+ \}_{M\in \MSmlog}\]
        where $\Sigma_P$ is defined as in \ref{sec:def-SigmaP-def-proj-is-bousfield}. Let $\Lambda$ be the set of pushout product maps $f\square g$ where $f\in \Lambda'$  and $g\in \{ (\partial \Delta^n)_+ \hookrightarrow (\Delta^n)_+\}$. To show that $\Hom_{\MMb}(\kX, \kZ)$  is ${\cI}^c$-projective fibrant, it's enough to show that for every generating cofibration $i =\id_M\wedge \iota_{n,+}\colon (h_M)_+\wedge (\partial\Delta^n)_+ \to (h_M)_+\wedge (\Delta^n)_+$, the push-out product of $i$ and any $f\in \Lambda$ is still a composition of pushouts of maps in $\Lambda$. For this, it is enough to notice that for every $M, M'$ and every $K, L \in \mathcal{S}_\bullet$, we have a canonical isomorphism
       \[ ((h_M)_+ \wedge K) \Daytensor_\bullet ((h_{M'})_+ \wedge L) \cong (h_M \Daytensor h_{M'})_+ \wedge (K\wedge L), \] 
 that replaces the fourth listed point in the proof of \cite[Lemma 2.18]{DRO}. To conclude, we have to show that for every $I^c$-weak equivalence $f\colon \kY\to \kY'$, every ${\cI}^c$-projective fibrant $\kZ$, and every cofibrant $\kX$,  the map $\Map(f\wedge \kX,\kZ  ) = \cS_\bullet( (f\wedge \kX)^c,\kZ )$ is a weak equivalence. By the above argument, $\Hom_{\MMb}(\kX, \kZ)$ is $I^c$-projective fibrant, so that the natural map between the simplicial function spaces
        \[\cS_\bullet(f^c\wedge \kX, \kZ)\to \cS_\bullet(f^c, \Hom_{\MMb}(\kX, \kZ) ) = \Map(f, \Hom_{\MMb}(\kX, \kZ)),\]
        induced by the closed monoidal structure on $\MMb$, is a weak equivalence. In particular, this shows that $f^c\wedge \kX$ is an ${\cI}^c$-projective weak equivalence. Since the ${\cI}^c$-local structure is obtained by left Bousfield localization from the projective (objectwise) model structure on simplicial presheaves, and the latter is monoidal with respect to the smash product, we conclude as in \cite{DRO} that $f^c\wedge \kX$  is an ${\cI}^c$-projective weak equivalence if and only if $(f\wedge \kX)^c$ is an ${\cI}^c$-projective weak equivalence. 
        
        For the general case, simply replace $\kX$ with $\kX^c\to \kX$, where $(-)^c$ denote a functorially chosen cofibrant replacement. The morphism $\kX^c\to \kX$ is an objectwise weak equivalence, so it is preserved by smashing with any motivic space. As for $\kX^c$ we can apply the previous claim.
        \end{proof} 
    \end{lem}
    \begin{prop}\label{prop:smash-preserves-we-cof-inj}The smash product preserves ${\cI}$-weak equivalences and cofibrations for the injective ${\cI}$-local model structure on $\MMb$, and induces a symmetric closed monoidal structure on the unstable motivic homotopy category $\MHp$.
        \begin{proof} By the description of weak equivalences in the pointed model category given by Theorem \ref{thm:pointed-version-model-cats} and the equivalence between the classes of ${\cI}^c$-projective weak equivalences and injective ${\cI}$-weak equivalences given by Proposition \ref{prop:proj-inj-local-agree}, Lemma \ref{lem:smash-preserves-I-proj-we} implies immediately that ${\cI}$-weak equivalences are preserved under smash product. Since cofibrations in the injective structure are monomorphisms, they are clearly preserved by smash product. This gives the homotopy category $\MHp$ the desired structure of monoidal category. To show that $\MHp$ is closed with respect to this monoidal structure, it's enough to use the fact that $\MMb^{{\cI}^c-loc}_{\rm proj}$ is a monoidal model category for the smash product. This is a consequence of Proposition \ref{prop:loc-monoidal-structure-CT}, using Lemma \ref{lem:smash-preserves-I-proj-we}, together with left properness of the ${\cI}^c$-local projective model structure on $\MMb$.
            \end{proof}
        \end{prop}
     \subsection{A  representability result} We conclude this Section with a general representability result in the ${\cI}$-homotopy category. 
     \subsubsection{} Let $\MHp$ be again the pointed unstable motivic homotopy category with modulus. For $\kX$ and $\kY$ motivic spaces with modulus, we set
     \[[\kX, \kY]_{\MHp} = \hom_{\MHp}(\kX, \kY),\]
(not to be confused with the internal hom for Day convolution, that we denoted $[-,-]$).     Let $M\in \MSmlog$ be any modulus datum. Evaluation at $M$ determines a Quillen pair 
     \begin{equation}\label{eq:ev-fr-adjunction}
         (M)_+\wedge (-) \colon \cS_\bullet  \leftrightarrows \MMb\colon {\rm Ev}_{M} = \cS_\bullet((M)_+, -)\end{equation}
  (where we write $M$ instead of $h_M$ for short) for the injective model structure $\MMb^{{\cI}-loc}_{\rm inj}$ on $\MMb$, since    $(M)_+\wedge (-)$ preserves monomorphisms and sectionwise weak equivalences.
    \subsubsection{} Write $S^1$ for the constant simplicial presheaf given by $\Delta^1/\partial \Delta^1$ and $S^n$ for the  $n$-th simplicial sphere $(S^1)^{\wedge n}$. The space $S^1$ is naturally pointed by the image of $\partial \Delta^1$, so that we can consider it as object in $\MMb$. Smashing with $S^1$ defines an endofunctor on $\MMb$, that is a left Quillen functor by Proposition \ref{prop:smash-preserves-we-cof-inj}. As customary, we write $\Sigma(-)$ for $S^1\wedge(-)$ (the suspension functor) and $\Omega^1(-) = \Hom_{\MMb}(S^1, -)$ for its right adjoint:
    \[\Sigma\colon \MMb\leftrightarrows \MMb \colon \Omega^1.\]
\begin{thm}\label{thm:formal-representability} Let $\kX$ be a pointed motivic space with modulus that satisfies the equivalent conditions of Proposition \ref{prop:BG}. Then for any pointed simplicial set $K$ and any modulus datum $M$, we have a natural isomorphism
    \[[K, \kX(M)]_{\cS_\bullet} \cong [(M)_+\wedge K,\kX]_{\MHp}.\]
    \begin{proof} The proof is a formal consequence of the results collected so far. Let $\varphi_\kX\colon \kX \to \cL \kX$ be a fibrant replacement for $\kX$ in the injective ${\cI}$-local model structure. By Proposition \ref{prop:BG}, $\varphi_\kX$ is a sectionwise weak equivalence. Then we have
        \[[K, \kX(M)]_{\cS_\bullet}\cong [K,  \cL \kX(M)]_{\cS_\bullet} \cong [K, \mathbf{R}{\rm Ev}_M ( \kX)]_{\cS_\bullet} \cong[(M)_+\wedge K,\kX]_{\MHp}  \]
        where the last isomorphism follows from the Quillen adjunction displayed in \eqref{eq:ev-fr-adjunction}.
        \end{proof}
    \end{thm}
    We state as separate Corollary the following result. It is a direct consequence of Theorem \ref{thm:formal-representability} and Corollary \ref{cor:boxlocal-implies-Ilocal}. This result shows that the category $\MHp$ can serve as an environment for studying representability problems for $\ol{\square}$-invariant (generalized) cohomology theories, as explained in the Introduction.
    \begin{cor}\label{cor:representability-box-case}Let $\kX$ be a pointed motivic space with modulus that is Nisnevich excisive in the sense of Definition \ref{def:Nisn-excisive-prop} and $\ol{\square}^1$-$\tensor$-invariant. Then for any $n\geq 0$ and any modulus datum $M$, we have a natural isomorphism
         \[\pi_n(\kX(M)) \cong [S^n\wedge (M)_+,\kX]_{\MHp}.\]
        \end{cor}
\noindent\emph{Acknowledgments.} This paper is essentially based on the second chapter of my PhD thesis, written in Essen under the supervision of Marc Levine. It is a pleasure to thank him heartily for much advice, constant support and encouragement. I would also like to thank Moritz Kerz and Shuji Saito for valuable comments and many friendly conversations on these topics, as well as Lorenzo Mantovani and Mauro Porta for many remarks on a preliminary version of this paper.  

Finally, I'm grateful to the anonymous referee for the careful reading of the manuscript, and for suggesting many improvements in the exposition. 
    \bibliography{bibHModulus} 

\begin{thebibliography}{10}

\bibitem{BThesis}
{\sc F.~Binda}, {\em Motives and algebraic cycles with moduli conditions}, PhD
  thesis, University of Duisburg-Essen, 2016.

\bibitem{BMilnorChowMod}
\leavevmode\vrule height 2pt depth -1.6pt width 23pt, {\em A cycle class map
  from {C}how groups with modulus to relative {$K$}-theory}, Doc. Math.,
  (2018), pp.~407--444.

\bibitem{BK}
{\sc F.~Binda and A.~Krishna}, {\em Zero cycles with modulus and zero cycles on
  singular varieties}, Compositio Math., 154 (2018), pp.~120---187.

\bibitem{BPO}
{\sc F.~Binda, D.~Park, and P.~A. {\O}stv{\ae}r}, {\em Triangulated categories
  of logarithmic motives over a field}, 2019.
\newblock In preparation.

\bibitem{BRS}
{\sc F.~Binda, K.~{R}{\"u}lling, and S.~Saito}, {\em On the cohomology of
  reciprocity sheaves},  (2019).
\newblock In preparation.

\bibitem{BS}
{\sc F.~Binda and S.~Saito}, {\em Relative cycles with moduli and regulator
  maps}, J.~of the Institute of Math.~Jussieu, to appear (2018), pp.~1--61.

\bibitem{Blander}
{\sc B.~A. Blander}, {\em Local projective model structures on simplicial
  presheaves}, $K$-Theory, 24 (2001), pp.~283--301.

\bibitem{BEAdditive}
{\sc S.~Bloch and H.~Esnault}, {\em The additive dilogarithm}, Doc. Math.,
  (2003), pp.~131--155 (electronic).
\newblock Kazuya Kato's fiftieth birthday.

\bibitem{BEAdditiveChow}
\leavevmode\vrule height 2pt depth -1.6pt width 23pt, {\em An additive version
  of higher {C}how groups}, Ann. Sci. {\'E}cole Norm. Sup. (4), 36 (2003),
  pp.~463--477.

\bibitem{CT}
{\sc D.-C. Cisinski and G.~Tabuada}, {\em Symmetric monoidal structure on
  non-commutative motives}, Journal of K-theory: K-theory and its Applications
  to Algebra, Geometry, and Topology, 9 (2012), pp.~201--268.

\bibitem{Day}
{\sc B.~Day}, {\em On closed categories of functors}, in Reports of the
  {M}idwest {C}ategory {S}eminar, {IV}, Lecture Notes in Mathematics, Vol. 137,
  Springer, Berlin, 1970, pp.~1--38.

\bibitem{DRO}
{\sc B.~I. Dundas, O.~R{\"o}ndigs, and P.~A. {\O}stv{\ae}r}, {\em Motivic
  functors}, Doc. Math., 8 (2003), pp.~489--525 (electronic).

\bibitem{GZFractions}
{\sc P.~Gabriel and M.~Zisman}, {\em Calculus of fractions and homotopy
  theory}, Ergebnisse der Mathematik und ihrer Grenzgebiete, Band 35,
  Springer-Verlag New York, Inc., New York, 1967.

\bibitem{Heller}
{\sc A.~Heller}, {\em Homotopy theories}, Mem. Amer. Math. Soc., 71 (1988),
  pp.~vi+78.

\bibitem{Hirsch}
{\sc P.~S. Hirschhorn}, {\em Model categories and their localizations}, vol.~99
  of Mathematical Surveys and Monographs, American Mathematical Society,
  Providence, RI, 2003.

\bibitem{Hovey}
{\sc M.~Hovey}, {\em Model categories}, vol.~63 of Mathematical Surveys and
  Monographs, American Mathematical Society, Providence, RI, 1999.

\bibitem{ImKelly}
{\sc G.~B. Im and G.~M. Kelly}, {\em A universal property of the convolution
  monoidal structure}, J. Pure Appl. Algebra, 43 (1986), pp.~75--88.

\bibitem{JardineJPAA}
{\sc J.~F. Jardine}, {\em Simplicial presheaves}, J. Pure Appl. Algebra, 47
  (1987), pp.~35--87.

\bibitem{JardineSpectra}
\leavevmode\vrule height 2pt depth -1.6pt width 23pt, {\em Motivic symmetric
  spectra}, Doc. Math., 5 (2000), pp.~445--553 (electronic).

\bibitem{JardineLocal}
\leavevmode\vrule height 2pt depth -1.6pt width 23pt, {\em Local homotopy
  theory}, Springer Monographs in Mathematics, Springer, New York, 2015.

\bibitem{KMSY}
{\sc B.~Kahn, H.~Miyazaki, S.~Saito, and T.~Yamazaki}, {\em Motives with
  modulus {I}: Modulus sheaves with transfers},  (2019).
\newblock Preprint.

\bibitem{KSY2}
{\sc B.~Kahn, S.~Saito, and T.~Yamazaki}, {\em Motives with modulus},  (2015).
\newblock arXiv:1511.07124v1.

\bibitem{KSY}
\leavevmode\vrule height 2pt depth -1.6pt width 23pt, {\em Reciprocity
  sheaves}, Compositio Mathematica, 152 (2016), pp.~1851--1898.
\newblock With two appendices by {K}ay {R}{\"u}lling.

\bibitem{KrishnaOnCycles}
{\sc A.~Krishna}, {\em On 0-cycles with modulus}, Algebra Number Theory, 9
  (2015), pp.~2397--2415.

\bibitem{KL}
{\sc A.~Krishna and M.~Levine}, {\em Additive higher {C}how groups of schemes},
  J. Reine Angew. Math., 619 (2008), pp.~75--140.

\bibitem{KP}
{\sc A.~Krishna and J.~Park}, {\em Moving lemma for additive higher {C}how
  groups}, Algebra Number Theory, 6 (2012), pp.~293--326.

\bibitem{KP4}
\leavevmode\vrule height 2pt depth -1.6pt width 23pt, {\em Algebraic cycles and
  crystalline cohomology},  (2015).
\newblock arXiv:1504.08181v4.

\bibitem{KP3}
\leavevmode\vrule height 2pt depth -1.6pt width 23pt, {\em A module structure
  and a vanishing theorem for cycles with modulus}, Math. Res. Lett., 24
  (2017), pp.~1147--1176.

\bibitem{LevineSmooth}
{\sc M.~Levine}, {\em Smooth motives}, in Motives and algebraic cycles, vol.~56
  of Fields Inst. Commun., Amer. Math. Soc., Providence, RI, 2009,
  pp.~175--231.

\bibitem{MV}
{\sc F.~Morel and V.~Voevodsky}, {\em {${\bf A}^1$}-homotopy theory of
  schemes}, Inst. Hautes \'Etudes Sci. Publ. Math.,  (1999), pp.~45--143
  (2001).

\bibitem{RulSaito}
{\sc K.~{R}{\"u}lling and S.~Saito}, {\em Higher {C}how groups with modulus and
  relative {M}ilnor {K}-theory}, Trans. Amer. Math. Soc.,  (2018),
  pp.~987--1043.
\newblock arXiv:1504.02669.

\bibitem{stacks-project}
{\sc {Stacks Project Authors}}, {\em {\itshape {S}tacks {P}roject}}.
\newblock \url{http://stacks.math.columbia.edu}, 2016.

\bibitem{TT}
{\sc R.~W. Thomason and T.~Trobaugh}, {\em Higher algebraic {$K$}-theory of
  schemes and of derived categories}, in The {G}rothendieck {F}estschrift,
  {V}ol.\ {III}, vol.~88 of Progr. Math., Birkh\"auser Boston, Boston, MA,
  1990, pp.~247--435.

\bibitem{VSelecta}
{\sc V.~Voevodsky}, {\em Homology of schemes}, Selecta Math. (N.S.), 2 (1996),
  pp.~111--153.

\bibitem{Vcdtop}
\leavevmode\vrule height 2pt depth -1.6pt width 23pt, {\em Homotopy theory of
  simplicial sheaves in completely decomposable topologies}, J. Pure Appl.
  Algebra, 214 (2010), pp.~1384--1398.

\bibitem{Vunst-nis-cdh}
\leavevmode\vrule height 2pt depth -1.6pt width 23pt, {\em Unstable motivic
  homotopy categories in {N}isnevich and cdh-topologies}, J. Pure Appl.
  Algebra, 214 (2010), pp.~1399--1406.

\end{thebibliography}
    \bibliographystyle{siam}
\end{document}